\documentclass[11pt]{article}
\usepackage{graphicx}
\usepackage{psfig}
\usepackage{latexsym}
\usepackage{amssymb}

\def\beq{\begin{eqnarray}}
\def\eeq{\end{eqnarray}}
\def\bqt{\begin{quotation}\vspace{-7pt}\noindent}
\def\eqt{\end{quotation}\vspace{-7pt}}
\def\Var{\mbox{Var}}

\begin{document}

\fontsize{11}{14.5pt}\selectfont

\begin{center}

{\small Technical Report No.\ 0607,
 Department of Statistics, University of Toronto}

\vspace*{0.45in}

{\LARGE\bf 
 Puzzles of Anthropic Reasoning Resolved \\[4pt]
 Using Full Non-indexical Conditioning\\[18pt]}

{\large Radford M. Neal}\\[3pt]
 Department of Statistics and Department of Computer Science \\
 University of Toronto, Toronto, Ontario, Canada \\
 \texttt{http://www.cs.utoronto.ca/$\sim$radford/} \\
 \texttt{radford@stat.utoronto.ca}\\[10pt]

 23 August 2006
\end{center}

\vspace{8pt} 

\noindent \textbf{Abstract.}  I consider the puzzles arising from four
interrelated problems involving ``anthropic'' reasoning, and in
particular the ``Self-Sampling Assumption'' (SSA) --- that one should
reason as if one were randomly chosen from the set of all observers in
a suitable reference class.  The problem of Freak Observers might
appear to force acceptance of SSA if any empirical evidence is to be
credited.  The Sleeping Beauty problem arguably shows that one should
also accept the ``Self-Indication Assumption'' (SIA) --- that one
should take one's own existence as evidence that the number of
observers is more likely to be large than small.  But this assumption
produces apparently absurd results in the Presumptuous Philosopher
problem.  Without SIA, however, a definitive refutation of the
counterintuitive Doomsday Argument seems difficult.  I\nolinebreak{}
show that these problems are satisfyingly resolved by applying the
principle that one should always condition on \textit{all} evidence
--- not just on the fact that you are an intelligent observer, or that
you are human, but on the fact that you are a human with a specific
set of memories.  This ``Full Non-indexical Conditioning'' (FNC)
approach usually produces the same results as assuming both SSA and
SIA, with a sufficiently broad reference class, while avoiding their
\textit{ad hoc} aspects.  I argue that the results of FNC are correct
using the device of hypothetical ``companion'' observers, whose
existence clarifies what principles of reasoning are valid.
I\nolinebreak{} conclude by discussing how one can use FNC to infer
how densely we should expect intelligent species to occur, and by
examining recent anthropic arguments in inflationary and string theory
cosmology.

\vspace*{8pt}

\section{\hspace*{-7pt}Introduction}\label{sec-intro}\vspace*{-10pt}

Accounting for selection effects is clearly necessary when drawing
conclusions from empirical data.  A poll conducted by telephone, for
example, will not tell us the opinions of people who don't have
telephones.  This simple observation has been seen by some, beginning
with Brandon Carter (1974), as having profound cosmological
implications, expressed as the Anthropic Principle --- ``what we can
expect to observe must be restricted by the conditions necessary for
our presence as observers''.  One typical cosmological application of
the Anthropic Principle is in ``explaining'' the observed values of
physical constants by assuming that they take on all possible values
in a multiplicity of universes, but that we of course must observe
values that are compatible with the existence of life.  A related use
of the Anthropic Principle is to deny that a cosmological theory in
which life is common should be considered more probable (other things
being equal) than one in which life is rare, as long as the latter
theory gives high probability to the existence of at least one
intelligent observer.

There is a large literature on the Anthropic Principle, much of it too
confused to address.  A coherent probabilistic account of the issues
involved has been presented by Nick Bostrom (2002, 2005), whose views
of the subject I consider worth critiquing.  Ken Olum's (2002, 2004)
views are also interesting, and are sometimes closer to my own.
Leonard Susskind (2006) and Leo Smolin (2006) have contrasting views
on the cosmological implications of the Anthropic Principle, which I
discuss at the end of this paper.

\subsection{\hspace*{-7pt}Four puzzles and two assumptions}\vspace*{-5pt}

One formalization of the intuition regarding
observer selection effects is what Bostrom calls the 
``Self-Sampling Assumption'' (SSA):
\bqt
  (SSA) One should reason as if one were a random sample from the
  set of all observers in one's reference class. (Bostrom 2002, p.~57)
\eqt
Bostrom regards this as a preliminary formulation; in particular,
he later considers more fine-grained ``observer moments''.  However, 
all forms of SSA require some specification of an appropriate 
``reference class'' (eg, all humans, or all intelligent observers), and hence 
are poorly defined if no precise basis for specifying such a class is given.

Despite this difficulty, something like SSA might appear to be essential
in order to deal with the Freak Observers problem:
\bqt
  How can vast-world cosmologies have \textit{any} observational
  consequences \textit{at all}?  We shall show that these cosmologies
  imply, or give very high probability to, the proposition that
  every possible observation is in fact made. (Bostrom 2002, p.~52)
\eqt
Bostrom argues that in a sufficiently large universe, brains in any
possible state will be emitted as Hawking radiation from black holes,
or condense from gas clouds as a result of large thermal fluctuations.
We need not consider such extreme possibilities in order to
see a problem, however.  Scientific experiments commonly have some
small probably of producing incorrect results for more mundane reasons.
In a large enough universe, it is likely that \textit{some} observer 
has made a misleading observation of any quantity of interest.
So, for example, that some observer in the universe
has made observations that with high confidence could be produced only
if the cosmic microwave background radiation is anisotropic is
\textit{no reason at all} to think that the background
radiation is actually anisotropic.  If we are to draw any conclusions
from observations we make, we need to see them not just as
observations that have been made, but as observations that have been
made \textit{by us}.  Bostrom argues that SSA together with the fact
that most observations made are not misleading then allows us to
conclude that our observations are likely to be correct.

However, if we accept SSA, we are led to the Doomsday Argument
expounded by John Leslie (1996), who attributes it to Carter.  The
Doomsday Argument says that your ordinary estimate of the chance of
early human extinction (based on factors such as your assessment of
the probability of an asteroid colliding with earth) should be
increased to account for an observer selection effect.  It is claimed
that the circumstance of your being (roughly) the 60 billionth human
to ever exist is more likely if there will never be more than a few
hundred billion humans than if there will be hundreds of trillions of
humans, as will be the case if humanity survives and colonizes the
galaxy.  This argument implicitly assumes SSA.  (If the reference
class is all intelligent observers, the argument requires that
uncertainty in time of extinction be shared with other intelligent
species.)

Although Leslie (1996), Carter (2004), and some others accept the
Doomsday Argument as valid, I take it to be absurd, primarily because
the answer it produces depends arbitrarily on the choice of reference
class.  Bostrom (2002) argues that this choice is analogous to a
choice of prior in Bayesian inference, which many are untroubled by.
However, a Bayesian prior reflects beliefs about the world.  A choice
of reference class has no connection to factual beliefs, but instead
reflects an ethical judgement, if it reflects anything.  It is thus
unreasonable for such a choice to influence our beliefs about facts of
the world.

The challenge is therefore to explain exactly why the Doomsday
Argument is invalid, without also destroying our ability to draw
conclusion from empirical data despite the possibility of freak
observers.  Many refutations of the Doomsday Argument have been
attempted, but as argued by Bostrom (2002, Chapter 7), most of these
refutations are themselves flawed.  In particular, it is not enough to
adduce plausible-sounding principles that if correct would defuse the
Doomsday Argument if these same principles produce unacceptable
results in other contexts.

One way of avoiding the conclusion of the Doomsday Argument is to
accept the ``Self-Indication Assumption'' (SIA) --- that we should
take our own existence as evidence that the number of observers in our
reference class is more likely to be large than small.  The effect of
SIA is to cancel the effect of SSA in the Doomsday Argument, leaving
our beliefs about human extinction unchanged from whatever they were
originally.  I refer to this combination as SSA+SIA, and to SSA with a
denial of SIA as SSA$-$SIA.  Bostrom (2002) argues that SIA cannot be
correct because of the Presumptuous Philosopher problem.  Consider two
cosmological theories, A and B, of equal plausibility in light of
ordinary evidence.  Suppose theory A predicts that there are about ten
trillion intelligent species in the universe, whereas theory B
predicts that there are only about ten intelligent species.  A
presumptuous philosopher who accepts SIA would decide that theory A
was a trillion times more likely than theory B, and would continue to
believe theory A despite virtually any experimental evidence against
it, since the chance that the experiments apparently refuting A were
fraudulently or incompetently performed, or produced misleading
results just by chance, is surely much greater than one in a trillion.

Denying SIA also seems as if it might lead to problems, however.  The
Sleeping Beauty problem (Elga 2000) sets up a situation in which the
flip of a coin determines whether an observer experiences a situation
once, if the coin lands Heads, or twice (the second time with no
memory of the first), if the coin lands Tails.  Logic analogous to
accepting SSA+SIA leads one to conclude that upon experiencing this
situation, the observer should believe with probability 1/3 that the
coin landed Heads, whereas SSA$-$SIA leads one to conclude that the
observer should assess the probability of Heads as being 1/2.
Although some have argued that 1/2 is the correct answer (Lewis 2001,
Bostrom 2006), the arguments that 1/3 is the correct answer appear to
me to be conclusive.  These include an argument based on betting
considerations, and another argument I detail below.  One might
therefore be reluctant to abandon SIA.

To summarize, accepting SSA+SIA produces answers regarding Freak
Observers, Sleeping Beauty, and the Doomsday Argument that I consider
reasonable, but seems to produce unreasonable results for the
Presumptuous Philosopher problem.  SSA$-$SIA also resolves the problem
of Freak Observers, and produces what might seem like reasonable
results for the Presumptuous Philosopher problem, but produces results
I consider wrong regarding Sleeping Beauty and the Doomsday Argument.

\subsection{\hspace*{-7pt}Resolving the puzzles}\vspace*{-5pt}

In this paper, I show that this dilemma can be resolved by abandoning
both SSA and SIA.  Both are \textit{ad hoc} devices with no convincing
rationale, and both require a ``reference class'' of observers, the
selection of which is quite arbitrary.  Instead, I advocate
consistently applying the general principle that one should condition
on \textit{all} the evidence available, including all the details of
one's memory, but without considering ``indexical'' information
regarding one's place in the universe (as opposed to what the universe
contains).  I call this approach ``Full Non-indexical Conditioning''
(FNC).

The results using FNC are the same as those found using SSA+SIA, when
it is clear how to apply the latter method.  As the problems I
consider will illustrate, however, FNC is a more general and more
natural method of inference, and has a clearer justification.

To test whether the conclusions found by using FNC or by using
alternative principles are correct, I introduce the device of
``companion'' observers.  For the Sleeping Beauty problem, this device
provides further evidence that the correct answer is obtained by FNC
(and by SSA+SIA), whereas the answer produced by applying SSA$-$SIA in
the manner previously done is incorrect.  Consideration of companion
observers also shows that SSA$-$SIA produces unacceptable results when
used with certain reference classes, including the narrow reference
classes that have previously been used for the Sleeping Beauty
problem.

When considering the Freak Observers and Presumptuous Philosopher
problems, I advocate restricting consideration to cosmological
theories in which the universe may be very large, but not so large
that it is likely to contain multiple observers with exactly the same
memories.  The problem of Freak Observers can then be resolved using
FNC, without any need for SSA.  I argue that as a general
methodological principle, one must be cautious of pushing thought
experiments to extremes, as this has produced spurious paradoxical
results in other contexts.  I do consider the possibility of infinite
universes later, in connection with inflationary cosmology.

I argue that there are actually two versions of the Presumptuous
Philosopher problem, with possibly different answers.  When comparing
theories differing in the density of observers, but not in the size of
the universe, consideration of companion observers provides good
reason to doubt the results found using SSA$-$SIA, whereas the results
of applying SSA+SIA or FNC appear correct.  I argue that no clear
conclusions can be drawn from the Presumptuous Philosopher problem
when the theories compared differ in the size of the universe.  The
Presumptuous philosopher problem therefore fails to provide a reason
to reject SSA+SIA or FNC.

\subsection{\hspace*{-7pt}Applying FNC to cosmology}\vspace*{-5pt}

After showing that FNC provides reasonable answers for each of the
four problems described above, I use FNC to estimate how densely we
should expect intelligent observers to occur in the galaxy.  This
discussion is not entirely \textit{a~priori}, but is based also on the
observed lack of extraterrestrials in our vicinity, both now, and as
far as we can tell, in the past.  The results I obtain shed some light
on the ``Fermi Paradox'' --- there are reasons to think
extraterrestrials should be common in the universe, but if so, where
are they?  My conclusions imply some pessimism regarding our future
prospects, but this is of a milder degree than that produced by the
Doomsday Argument, and follows from empirical evidence, not anthropic
reasoning.

I conclude by discussing the implications of FNC for anthropic
arguments relating to inflationary cosmology, which favours a universe
or universes of infinite extent, and to cosmologies based on string
theory, in which a multiplicity of universes populate a ``landscape''
of differing physical laws.

\section{\hspace*{-7pt}Methodology}\label{sec-method}\vspace*{-10pt}

Before discussing the four problems of anthropic reasoning mentioned
above, a general examination of the methodology to be employed seems
desirable.

\subsection{\hspace*{-7pt}The nature of probabilities}\vspace*{-5pt}

First, since all these problems involve probabilistic answers, one may
ask what these probabilities mean.  I interpret probabilities as
justified subjective degrees of belief --- subjective in that they
depend on the information (including prior information) available to
the subject, and justified in that they follow from correct principles
of reasoning, rather than being capricious.  

Probabilistic beliefs about scientific hypotheses (eg, whether or not
earth-like planets are common in the universe) are based partly on our
prior assessments of plausibility.  Often, such hypotheses are guesses
about the implications of some more fundamental theory, whose true
implications cannot be computed exactly, but for which approximations
and mathematical intuition provide some guide.  We modify these prior
beliefs according to how successfully these hypotheses account for
observations (eg, of whether nearby stars have planets).

The probabilistic nature of most predictions may be due to at least four 
sources:\vspace{-6pt}
\begin{enumerate}
\item[1)] Inherent randomness in physical phenomena.
\item[2)] Ignorance about the initial conditions of physical phenomena.
\item[3)] Ignorance about our place in the universe.
\item[4)] Inability to fully deduce the consequences of a theory.\vspace*{-6pt}
\end{enumerate}
These possible sources of uncertainty are not mutually exclusive.  
At least (1), (2), and (4) are common sources of uncertainty in
ordinary scientific reasoning.  One's interpretation of
quantum mechanics determines whether ``random'' quantum phenomena
are seen as examples of (1), in the Copenhagen interpretation, 
or of (3), in the Many Worlds interpretation.  Leaving aside the
technicalities of the interpretation of this particular theory, one 
might generally take the ontological position that any apparently random 
choice is actually made in all possible ways, in parallel universes, all
of which are real (though perhaps with different ``weights'', 
corresponding to the probabilities of the choices), thereby
converting physical randomness to ignorance about which parallel
universe we are in.  

In trying to resolve the puzzles of anthropic reasoning addressed in
the paper, it seems best to not also attempt to resolve issues
regarding the nature of probability in physical theories.  I will look
for a solution based on fairly common sense notions, presuming that
these will in essence survive any final resolution of issues such as
the interpretation of quantum mechanics.  In thought experiments, I
will follow convention by usually talking about choices that are
determined by a coin flip, whose randomness likely derives from
source~(2).  The reader may, however, replace this with uncertainty of
another type, such as whether the 1,341,735'th digit of $\pi$ is even
or odd, assuming that this digit is not already known to the people
involved.

\subsection{\hspace*{-7pt}Indexical information and 
            reference classes}\vspace*{-5pt}\label{sec-idx}

A central feature of anthropic reasoning is the use of ``indexical''
information that certain observations were not just made, but were
made \textit{by you}.  Another expression of this concept is that you
should consider not just ``possible worlds'', but ``possible
\textit{centred} worlds'', in which your location in the universe is
specified (Lewis 1979).

Both SSA and SIA involve the indexical information that \textit{you}
are members of a certain reference class.  The conclusions that follow
from these principles individually are sensitive to the choice of this
reference class.  Interestingly, however, when \textit{both} SSA and
SIA are assumed, this dependence disappears, as long as the reference
class is broad enough to encompass all observers who could possibly
have observed what you have observed.

Let $C$ and $C'$ be two references classes of observers.  Let $D$ be
the set of observers who have observed the same data as you have.  I
will assume here that $D \subseteq C$ and $D \subseteq C'$ --- ie, the
data you observed indicates that you yourself are a member of both of
these reference classes.  I will also use the symbol $D$ to denote the
event that you are in the set $D$.  Let $A$ and $B$ be two
mutually-exclusive hypotheses about the universe, each of which 
specifies the numbers of observers in $C$, $C'$, and $D$.  Suppose that
based on prior information of the usual sort, $P(A) = P(B) = 1/2$.  If
$A$ is true, the number of observers in class $C$ is $|C|_A$ and the
number in class $C'$ is $|C'|_A$; if $B$ is true, these numbers are
$|C|_B$ and $|C'|_B$.  If $A$ is true, the number of observers in $D$ is
$|D|_A$; if $B$ is true, this number is $|D|_B$.

If we assume SSA with reference class $C$, but not SIA, the
probability of hypothesis $A$ given the observed data is\vspace*{-4pt}
\beq
  P(A\,|\,D) & = & { P(A)\,P(D\,|\,A) \over 
                     P(A)\,P(D\,|\,A) \ +\ P(B)\,P(D\,|\,B)} \\[6pt]
             & = & { (1/2)\, (|D|_A/|C|_A) \over
                     (1/2)\, (|D|_A/|C|_A)\  +\ (1/2)\, (|D|_B/|C|_B) } 
\\[-6pt]\nonumber
\eeq
If we instead use reference class $C'$, $P(A\,|\,D)$ will be given
by this formula with $|C'|_A$ and $|C'|_B$ replacing $|C|_A$ and $|C|_B$.
The result will in general be different.

However, if we assume SSA+SIA, the prior probabilities for 
$A$ and $B$ are modified in proportion to the number of observers
they imply are in the reference class.  This gives the following
formula for $P(A\,|\,D)$ when the reference class is $C$:
\beq
  P(A\,|\,D) &\!\! = \!\!& { (|C|_A/2)\, (|D|_A/|C|_A) \over
                     (|C|_A/2)\, (|D|_A/|C|_A)\  +\ (|C|_B/2)\, (|D|_B/|C|_B) }
             \ = \ { (1/2)\, |D|_A \over (1/2)\,|D|_A\ +\ (1/2)\,|D|_B }
\ \ \ \ \ \ 
\eeq
Since $|C|_A$ and $|C|_B$ cancel, the identical result is obtained if $C'$ is 
used as the reference class.

This lack of dependence on the reference class suggests that even if
the right result is obtained by assuming both SSA and SIA, the joint
affirmation of these two principles may not be the most illuminating
way of describing the logic leading to this result.

\subsection{\hspace*{-7pt}Full Non-indexical Conditioning 
  (FNC)}\label{sec-fnc}\vspace*{-5pt}

I advocate probabilistic reasoning by the standard method of fully
conditioning on all information that you possess.  Of course, in most
ordinary circumstances, you can ignore much information that you know
is not relevant to the problem --- eg, when predicting tomorrow's
weather, you should condition on the current barometric pressure (if
you know it), but there is no need to also condition on the name of
your kindergarten teacher (even if you still remember it).  When in
doubt, however, it is always correct to conditional on additional
information, since if this information is in fact irrelevant,
conditioning on it will not change the result.

When dealing with puzzling instances of anthropic reasoning, what is
relevant and irrelevant is unclear, so I maintain that you should
condition on \textit{everything} you know --- your entire set of
memories --- to be sure of getting the right answer.  You might also
condition on the indexical information that these are \textit{your}
memories.  However, I will here consider what happens if you ignore
such indexical information, conditioning only on the fact that
\textit{someone} in the universe with your memories exists.  I refer
to this procedure as ``Full Non-indexical Conditioning'' (FNC).

Some ordinary situations might appear to require use of indexical
information, but a closer examination shows that FNC produces the
correct answer in these cases.  For instance, suppose that you and
some number of other people are recruited as subjects for an
experiment.  You do not know the number of subjects for this
experiment, but based on your knowledge of the budget limitations in
the field of experimental philosophy, you have a prior distribution
for this number, $N$, that is uniform over the integers from 1 to 20.
You and the other subjects are taken to separate rooms, without seeing
each other, where you are instructed to flip a fair coin three times and
record the sequence of Heads and Tails obtained.  You record the
sequence HTT.  On the basis of this new data, what should be your
posterior distribution for $N$? 

An invalid way of reasoning here is to condition on the fact that the
sequence HTT was recorded by a subject of the experiment, and
on that basis conclude that your posterior distribution should 
be
\beq
P(N=n \,|\, \mbox{HTT recorded})
   & = &
   {P(N=n)\ P(\mbox{HTT recorded}\,|\,N=n) \over
   \sum_{n'=1}^{\infty} P(N=n')\ P(\mbox{HTT recorded}\,|\,N=n')} \\[5pt]
   & = &
   {P(N=n)\ (1 \ -\ P(\mbox{HTT not recorded}\,|\,N=n)) \over
   \sum_{n'=1}^{\infty} P(N=n')\ 
           (1 \ -\ P(\mbox{HTT not recorded}\,|\,N=n'))}\ \ \ \ \ \\[5pt]
   & = &
   {(1/20)\ (1 \ -\ (1 - 2^{-3})^n) \over
    \sum_{n'=1}^{20}\, (1/20)\ (1 \ -\ (1 - 2^{-3})^{n'})} \\[-6pt]\nonumber
\eeq
These probabilities vary from 0.0093 for $n=1$ to 0.0690 for $n=20$.
Intuitively, these probabilities seem wrong, for two reasons.  First,
since you had to record \textit{some} sequence of flips, it seems that
knowledge of the particular sequence you recorded shouldn't change your 
beliefs about $N$.
On the other hand, it seems that the fact that \textit{you} were recruited for 
the experiment should increase the probability that there are many subjects
(by more than happens above).

The problem is fixed if you take account of the indexical information
that it was \textit{you} who recorded the sequence \mbox{HTT}.
The probability of this sequence being obtained by you is $2^{-3}$
regardless of $N$.  On the other hand, if the pool of possible
experimental subjects is of size $M$ (assumed to be greater
than 20), the probability that you will be recruited as a subject if
there are $n$ subjects is $n/M$.  The posterior distribution for $N$ that 
results is
\beq
P(N=n\,|\,\mbox{You were recruited and recorded HTT})
  & = & {(1/20)\ (n/M)\ 2^{-3} \over
         \sum_{n'=1}^{20}\, (1/20)\ (n'/M)\ 2^{-3}} \\[5pt]
  & = & \textstyle n\ /\ \sum\nolimits_{n'=1}^{20}\, n' \ \ =\ \ n\ /\ 210
\eeq
These probabilities vary from 0.0048 for $n=1$ to 0.0952 for $n=20$.

This answer is also obtained if we condition on \textit{all}
non-indexical information.  We know not just that the sequence HTT
was recorded, but also that it was recorded by a subject of your age,
with your hair colour, who went to a school just like yours, who has
your taste in music, who has the same opinion of rice pudding as you
do, etc.  The probability that the $i$'th subject recruited for the
experiment will have all these characteristics is some very small
number, $\epsilon$.  The probability that the $i$'th subject has these
characteristics and also records coin flips of HTT is 
$\epsilon\, 2^{-3}$.
Since this probability is extremely small, the probability that \textit{any}
of $n$ subjects will have these characteristics and record those flips
is very well approximated by $n\,\epsilon\,2^{-3}$.  Conditioning on
all non-indexical information therefore produces the following
posterior distribution for $N$:
\beq
P(N=n\,|\,\mbox{all non-indexical information})
  & = & {(1/20)\ n\, \epsilon\, 2^{-3} \over
         \sum_{n'=1}^{20}\, (1/20)\ n\, \epsilon\, 2^{-3}} \\[5pt]
  & = & \textstyle n\ /\ \sum\nolimits_{n'=1}^{20}\, n' \ \ =\ \ n\ /\ 210
\eeq
This is the same result as found above using indexical information.

Use of indexical information therefore seems unnecessary in ordinary
situations, since the non-indexical information regarding your
memories is normally sufficient to uniquely identify you.  SSA can
perhaps be seen as arising from what might be called ``Full Indexical
Conditioning'', in which we assume that in addition to their memories,
everyone also has some unique ``essence'', and everyone in some sense
knows what their own essence is.  One could then argue that you should
condition not only on your memories, but also on having your own
essence, which by assumption is shared with no one else.  Conditioning
on more than FNC rather than less seems the only way of reconciling
SSA$-$SIA with the fundamental principle that probabilistic
inferences should be based on all known information.  I will
illustrate this idea when discussing Sleeping Beauty in
Section~\ref{sec-beauty-fnc}.  It seems preferable to me to not
introduce such mystical ``essences'' unless ignoring them can be shown
to produce implausible results.  Moreover, if one does think in terms
of such essences, it seems hard not to proceed to acceptance of SIA as
well as SSA, which effectively renders thoughts of essences pointless,
as I will now explain.

It turns out that if the universe is not excessively huge, using
SSA+SIA produces the same results as using FNC.  Consider the
situation with two hypotheses, $A$ and $B$, discussed above in
Section~\ref{sec-idx}.  Assume now that you condition on all
information you remember, and that these memories are extensive enough
that there is only a small probability that an observer with your
memories would exist, under either hypothesis.  Let $|C|_A$ and
$|C|_B$ be the numbers of observers in some suitable reference class
(whose members ``might have had'' your memories) if hypotheses $A$ and
$B$ are true, respectively.  Let $\epsilon_A$ and $\epsilon_B$ be
the (extremely small) probabilities that a particular observer in this
reference class will have your memories under hypotheses $A$ and $B$.
Suppose that you assess the prior probability of $A$ and $B$ as
$P(A)=P(B)=1/2$, where this prior is based on your scientific
knowledge, but not on the multitude of details of your life that make
you unique.  Applying SIA will shift these priors to
$P(A)\,=\,(|C|_A/2)\,/\,(|C|_A/2+|C|_B/2)$ and
$P(B)\,=\,(|C|_B/2)\,/\,(|C|_A/2+|C|_B/2)$.  Applying SSA, the
probability that you will have your memories is $\epsilon_A$ if
hypothesis $A$ is true, and $\epsilon_B$ if hypothesis $B$ is true.
The result of applying both SIA and SSA is therefore
\beq
 \lefteqn{P(A\,|\,\mbox{all your memories})}\ \ \ \ \ \ \ \ \ \ \ \nonumber 
 \\[2pt]
  & = &
    { \epsilon_A\ (|C|_A/2)\,/\,(|C|_A/2+|C|_B/2) \over
      \epsilon_A\ (|C|_A/2)\,/\,(|C|_A/2+|C|_B/2)\ +\ 
      \epsilon_B\ (|C|_B/2)\,/\,(|C|_A/2+|C|_B/2)}\ \ \ \\[5pt]
  & = & { \epsilon_A\, |C|_A \over \epsilon_A\, |C|_A \ +\ \epsilon_B\, |C|_B}
  \\[-10pt] \nonumber
\eeq
Provided that $\epsilon_A\,|C|_A$ and $\epsilon_B\,|C|_B$ are both close to
zero (as they will be if the universe is not excessively huge), the same result 
will be obtained by applying FNC --- $\epsilon_A\,|C|_A$
and $\epsilon_B\,|C|_B$ are then very close to the probabilities of \textit{any}
observer with your memories existing under hypotheses $A$ and $B$, respectively,
which given the equal prior probabilities of $A$ and $B$ produces the same
result as above found using SSA and SIA. 

FNC is a more general principle of inference than SSA and SIA,
however, since it does not require any notion of a reference class.
FNC requires only that there be some way of computing the probability
that an observer with your memories will exist.  As was done above, it
is convenient to separate your scientific memories (which may be
shared with many others) from the rest of your memories, which make
you a unique observer.  Conditioning on your scientific memories
converts whatever primitive prior distribution you had regarding
scientific theories to what would ordinarily be regarded as your
prior.  We can then consider how this prior is altered by conditioning
on subsequent scientific observations and on the memories that make
you unique.

Note that the probabilities involved in FNC need not derive from some
random physical process, but may simply reflect ignorance or an
inability to fully deduce the consequences of known facts.  This will
be discussed further in Section~\ref{sec-density-fnc}.

\subsection{\hspace*{-7pt}Assessing arguments by considering companion 
                          observers}\label{sec-companions}\vspace*{-5pt}

If your opinions differ from those of an intelligent friend who
possesses the same information as you, you should question the
validity of the reasoning that led you to these opinions.  Ultimately,
after exchanging information and fully discussing the matter with your
friend, you should expect to come to the same conclusions regarding
factual matters.  Persistent disagreements might seem possible due to
differing prior beliefs, but as discussed by Hanson (2006), this is
possible for fully rational observers only if they disagree about the
processes by which they came to hold these prior beliefs.  Agreement
with hypothetical friends has been used as a test of reasoning in the
past --- in particular, Nozick (1969) uses it in his discussion of
Newcomb's Problem, as described in the next section.

I propose here to test the validity of anthropic arguments by
comparing the conclusions of such arguments with those that imaginary
``companion'' observers would reach, using the same principles of
reasoning (eg, acceptance of SSA but not SIA).  Of course, it is
possible that several incompatible sets of principles might each lead
to consistency with the conclusions of companions reasoning by these
same principles, so this test may sometimes fail to fully resolve the
issues.

By considering possible companions, a general constraint on the use of
SSA$-$SIA can be derived.  Suppose there exist two types of
observers, $X$ and $Y$.  These observers are considering two theories,
$A$ and $B$, according to which the numbers of observers of these
types are $|X|_A$ and $|Y|_A$, for theory $A$, and $|X|_B$ and
$|Y|_B$, for theory $B$.  There is a pairing observers of type $X$ with
companion observers of type $Y$, in which each observer is paired with at
most one companion observer.  If $|X|=|Y|$, all observers have
companions; otherwise, some observers of the more numerous type are
unpaired.

Suppose that all observers consider theories $A$ and $B$ to be equally
likely based on the usual sorts of evidence (ie, without applying SSA
or SIA).  Now consider what those observers with companions will
conclude by applying SSA$-$SIA using as their reference class only
observers of their own type, and taking account of their knowledge
that they were paired with a companion observer of the other type.
Observers of type $X$ will reason that their chance of having a
companion observer is $\min(1,|Y|_A/|X|_A)$ if theory $A$ is true,
and $\min(1,|Y|_B/|X|_B)$ if theory $B$ is true.  The odds an observer
of type $X$ assigns to theory $A$ over theory $B$ will therefore be
\beq
   {\min(1,|Y|_A/|X|_A) \over \min(1,|Y|_B/|X|_B)} & = &
   \left\{\begin{array}{ll}
     1               & \mbox{if $|Y|_A \ge |X|_A$ and $|Y|_B \ge |X|_B$} \\[4pt]
     |Y|_A\ /\ |X|_A & \mbox{if $|Y|_A \le |X|_A$ and $|Y|_B \ge |X|_B$} \\[4pt]
     |X|_B\ /\ |Y|_B & \mbox{if $|Y|_A \ge |X|_A$ and $|Y|_B \le |X|_B$} \\[4pt]
     |X|_B|Y|_A\ /\ |X|_A|Y|_B 
        & \mbox{if $|Y|_A \le |X|_A$ and $|Y|_B \le |X|_B$} 
   \end{array}\right.\ \ \ \
\label{eq-comp-odds-X}
\eeq
whereas for an observer of type $Y$, the odds in favour of theory $A$
would be
\beq
   {\min(1,|X|_A/|Y|_A) \over \min(1,|X|_B/|Y|_B)} & = &
   \left\{\begin{array}{ll}
     |X|_A|Y|_B\ /\ |X|_B|Y|_A                
        & \mbox{if $|Y|_A \ge |X|_A$ and $|Y|_B \ge |X|_B$} \\[4pt]
     |Y|_B\ /\ |X|_B & \mbox{if $|Y|_A \le |X|_A$ and $|Y|_B \ge |X|_B$} \\[4pt]
     |X|_A\ /\ |Y|_A & \mbox{if $|Y|_A \ge |X|_A$ and $|Y|_B \le |X|_B$} \\[4pt]
     1
        & \mbox{if $|Y|_A \le |X|_A$ and $|Y|_B \le |X|_B$} 
   \end{array}\right.\ \ \ \
\label{eq-comp-odds-Y}
\eeq
Since these are generally not equal, we see that companions will disagree
in this scenario if they each reason with SSA$-$SIA using as their
reference class only observers of their own type.  

However, these companion observers will agree if they apply SSA$-$SIA
using as their reference class all observers of both types, because of
the effects of a Doomsday-like argument, of a sort discussed further
in Section~\ref{sec-gendoom}.  Observers of type $X$ will reason that
their chances, applying SSA, of being of type $X$ are
$|X|_A/(|X|_A+|Y|_A)$ if theory $A$ is true, and $|X|_B/(|X|_B+|Y|_B)$
if theory $B$ is true.  The odds in favour of theory $A$ are therefore
multiplied by $(|X|_A/|X|_B)\ \times\ (|X|_B+|Y|_B)/(|X|_A+|Y|_A)$.  Multiplying
equation~(\ref{eq-comp-odds-X}) by this factor, we get that an observer
of type $X$ having a companion will consider the odds in favour of theory $A$
to be
\beq
   {|X|_B+|Y|_B \over |X|_A+|Y|_A}\ \times\
   \left\{\begin{array}{ll}
   |X|_A\ /\ |X|_B\ & \mbox{if $|Y|_A \ge |X|_A$ and $|Y|_B \ge |X|_B$} \\[4pt]
   |Y|_A\ /\ |X|_B\ & \mbox{if $|Y|_A \le |X|_A$ and $|Y|_B \ge |X|_B$} \\[4pt]
   |X|_A\ /\ |Y|_B\ & \mbox{if $|Y|_A \ge |X|_A$ and $|Y|_B \le |X|_B$} \\[4pt]
   |Y|_A\ /\ |Y|_B\ & \mbox{if $|Y|_A \le |X|_A$ and $|Y|_B \le |X|_B$} 
   \end{array}\right\}\ \ \ \
\eeq
Similarly, multiplying equation~(\ref{eq-comp-odds-Y}) by
$(|Y|_A/|Y|_B)\ \times\ (|X|_B+|Y|_B)/(|X|_A+|Y|_A)$ gives the 
odds in favour of theory $A$ for an observer of type $Y$ with a companion,
which turn out to be identical to the odds above.  

This computation shows that requiring consistency with conclusions of
a companion imposes a constraint on the reference class used with
SSA$-$SIA --- the companions must use the same reference class,
which must therefore include both of them.  This constraint might be
seen as making anthropic arguments based on SSA$-$SIA less arbitrary,
and hence more attractive.  However, this constraint also makes it
harder to apply such arguments, since to chose a suitable reference
class, you must know the full set of observers with whom you
would expect to agree.

In contrast, SSA+SIA produces consistent results even when observers
use reference classes that include only their own type of observer,
excluding their companion.  This may be confirmed by multiplying
equation~(\ref{eq-comp-odds-X}) by $|X|_A/|X|_B$ and
equation~(\ref{eq-comp-odds-Y}) by $|Y|_A/|Y|_B$, the factors by which
SIA modifies the prior odds.  Companion observers applying FNC will
obviously produce consistent conclusions, since FNC does not involve
indexical information, and companions are assumed to share all
non-indexical information.

\subsection{\hspace*{-7pt}The dangers of fantastic 
   assumptions}\label{sec-fant}\vspace*{-5pt}

Several of the puzzles treated here employ thought experiments, and
make other arguments, that are based on hypothetical and perhaps
fantastic assumptions.  This can sometimes produce spurious
conclusions.  We may accept a fantastic assumption, on the basis that
although it isn't true in reality, it ``might be true'', and then
proceed to reason utilizing other premises that are based on the
reality that we have implicitly rejected in making the fantastic
assumption.

Searle's (1980) Chinese Room Argument provides one example.  He argues
that a computer cannot possibly understand Chinese.  Any program that
enabled a computer to understand Chinese could in principle be
executed by a person in a room who takes inputs from a window, follows
certain simple rules for shuffling tokens about, and shoves results
out another window.  The person executing this program need have no
understanding of Chinese to begin with, and is unlikely to acquire any
understanding of Chinese by performing the tasks needed to execute
this program.  So, the argument goes, a computer running such a
program will also not really understand Chinese, regardless of whether
it might superficially appear to.

A common (and in my view, correct) response is that although the person
executing the program does not understand Chinese, the system of
person plus room is a physical embodiment of another entity that
\textit{does} understand Chinese.  To this, a defender of the Chinese
Room Argument may reply that, in principle, the room is unnecessary ---
a person with a sufficiently good memory could execute the program
entirely in their head, without the need of any physical tokens.  To
the subsequent objection that this just means that the new entity is
physically contained in the same body as the original person, Harnad
(2001) has mockingly replied 
\bqt
\noindent This was tantamount to conjecturing that, as a result of memorizing
and manipulating very many meaningless symbols, Chinese-understanding
would be induced either consciously in Searle, or,
multiple-personality-style, in another, conscious
Chinese-understanding entity inside his head of which Searle was
unaware.\vspace{5pt}

\noindent I will not dwell on any of these heroics; suffice it to say that even
Creationism could be saved by ad hoc speculations of this order.
\eqt
In a seminar I attended in the 1990's, Harnad explained in more
detail that the psychiatric literature on multiple personality 
disorder contains no recorded case of such a second being, with 
totally different language and other capabilities, existing within
someone's head.

This is an extreme case of making a fantastic assumption and then
reasoning with premises that contradict it.  Any program capable of
appearing to understand Chinese will very likely require a computer at
least as powerful as those available today to execute in real time.
Compared to manual execution by a person, today's computers are at
least a billion times faster, and have at least a billion times
as much readily-available memory.\footnote{I refer here to the speed
and memory available when a person consciously carries out the simple
tasks needed to execute the steps of a computer program.  The
computational power that underlies unconscious functions of our brains
likely exceeds that of today's computers.}  The characteristics of a
hypothetical person whose computational abilities are a billion times
greater than those of ordinary people are certainly not going to be
obvious to us, or deducible from the current psychiatric literature.
To assume the existence of such a person and then claim incredulity at
a consequence of their existence abuses the hospitality of one's
interlocutor in conceding that ``in~principle'' a computer program can
be executed manually --- when in reality, this is true only of
programs no more than a few pages long, that operate for no more than
a few hundred steps.

A more subtle example of the dangers of fantastic thought experiments
is provided by Newcomb's Problem, first discussed in print by Nozick
(1969).  We imagine that a wealthy ``Predictor'', who is very good at
predicting human behaviour, conducts a ``game'' that operates as
follows.  A person is randomly selected to participate, and is then
shown two boxes.  They are told that the first box contains \$1000,
and the other box contains either \$1,000,000 or \$0.  The participant
may either take both boxes, or take just the second box, and receives
all the money in the box or boxes they take.  The Predictor puts
\$1,000,000 in the second box if and only if he predicts that the
participant will take only this box.  Suppose that you have seen the
game played many times, and are convinced that the Predictor's
predictions are almost certain to be correct.  If you are selected to
be a participant in this game, should you take both boxes, or only the
second box?

The argument for taking only the second box is that you are then
almost certain to receive \$1,000,000, whereas if you take both boxes,
you almost certainly will receive only \$1000.  (Conceivably, you
might receive \$1,001,000 if you take both boxes, but only if the
Predictor is wrong, which you know is very unlikely.)  The argument
for taking both boxes is that you will then receive \$1000 more than
you would if you took only the second box, regardless of whether the
Predictor was right or wrong in his prediction, which he has
\textit{already} made.  Nozick (1969) favours taking both boxes,
strengthening this argument by pointing out that if a friend of yours
could see into both boxes, they would certainly advise you to take
both of them.  Actually, your friend cannot see into the boxes, or if
they can, they aren't allowed to advise you, but you know what advice
your friend would give if they could, and you should follow this
advice.

I would find the argument for taking two boxes convincing, if it were
not for a matter that seems to have been overlooked in the
philosophical literature\footnote{I wrote a paper on this idea twenty years
ago, which was rejected by \textit{Mind}, though verbal discussions
with philosophers over the years have been more positive.  At least
two other people have thought of this idea independently --- Scott
Aaronson described the idea in a November 2005 blog posting at
\texttt{http://www.scottaaronson.com/blog/2005/11/dude-its-like-you-read-my-mind.html},
and someone going by the name of ``Count Iblis'' wrote about it
in a December 2005 blog posting at \texttt{http://countiblis.blogspot.com}.} 
--- How does the Predictor make such accurate predictions?

Superficially, assuming the existence of such an accurate Predictor
may not seem too extreme.  We all predict other people's actions every
day, often successfully.  However, we also have a strong sense that we
have free will, and that our will is integrated with our whole being
--- for example, that any part of our memories can potentially affect
our actions.  Hence, while some of our actions are easy to predict,
other actions could only be predicted with high accuracy by a being
who has knowledge of almost every detail of our memories and
inclinations, and who uses this knowledge to simulate how we will act
in a given situation.  To predict whether a participant will take one
box or two (which for at least some people must be a difficult
behaviour to predict), the Predictor must have some way of measuring
with high accuracy the relevant aspects of the participant's brain (at
some time prior to when the game is played) and a very powerful
computer that can simulate the participant's mental
processes.\footnote{Note that we needn't assume that this simulation
is absolutely accurate, if we allow that the Predictor may be wrong
with some tiny probability.  Quite large final error rates would be
compatible with the Newcomb scenario as long most errors are
introduced by the Predictor's careless assistant, after a highly
accurate simulation has been run.  Such errors would be uncorrelated
with the type of participant, and hence you would not be justified in
feeling that you in particular might be able to ``beat the game''.
However, errors due to faulty simulation that arise when, for example,
a participant thinks of one particular argument would undermine the
Newcomb scenario.}

Now we can see why Newcomb's Problem involves an extreme fantastic
assumption --- the only plausible mechanism for accurate prediction
involves brain measurements and simulations that are far beyond our
current ability, and that may be impossible in principle, if quantum
effects are crucial to how the brain works (since non-destructive
copying of quantum states is not possible).\footnote{One might object
that some other mechanism for accurate prediction not involving
simulation might be possible.  However, Nozick specifically excludes
mechanisms, such as time travel, that introduce backward causation.
The onus is on someone wishing us to take the problem seriously as a
paradox to provide at least a hint of how such accurate predictions
might be obtained with neither accurate simulation nor backward
causation.}  This fantastic assumption has a crucial consequence ---
the simulation the Predictor conducts in order to predict your choice
will (if you accept a functionalist view of consciousness) create
another conscious being, and \textit{you have no way of knowing that
you are not this being}.  If you are the being in the simulation, your
``choice'' has a causal effect on whether the Predictor puts
\$1,000,000 or \$0 in the second box.  Supposing that the simulated
``you'' has sympathy for the real ``you'', or perhaps that ``you''
intended to donate the money to a worthy charity all along, it is now
clear that you should take only the second box, since that may cause
the real ``you'' to obtain \$1,000,000, and at worse costs the real
``you'' only \$1000.  Note that the argument involving advice from a
friend loses its force once the situation is really understood.  Your
friend may not actually be there (if you are being simulated, but he
is not), and if he is (and understands the situation), he will advise
you to take only the second box.

Possible problems with overly fantastic assumptions arise with several
of the problems discussed below.  We should also be careful to keep
the ``companion'' observers of the previous section from becoming too
fantastic, at least with respect to their cognitive and other relevant
characteristics (though other fantastic aspects may be innocuous).

\section{\hspace*{-7pt}Sleeping Beauty}\label{sec-sb}\vspace*{-10pt}

The Sleeping Beauty problem is described by Elga (2000).  On Sunday,
Beauty is put to sleep.  On Monday she is woken, then later put to
sleep again.  While she is awake, she does not have access to any
information that would help her infer the day of the week.  If a flip
of a fair coin lands Tails, Beauty is woken again on Tuesday, but only
after she is administered a drug that causes her to forget her Monday
awakening (leaving her memories in the same state as they were after
falling asleep on Sunday).  Again, she obtains no information that
would reveal the day of the week.  Regardless of how the coin lands,
Beauty is woken on Wednesday, and immediately told that the experiment
is over.  Beauty knows that this is how the experiment is set up.
When Beauty wakens before Wednesday, what probability should she
assign to the coin landing Heads?

\subsection{\hspace*{-7pt}The $1/2$ and $1/3$ answers}\vspace*{-5pt}

Some (eg, Lewis 2001) argue that on Sunday Beauty should certainly
assign probability 1/2 to the coin landing Heads, that upon wakening
she acquires no additional information (since she knew that she would
experience such an awakening regardless of how the coin lands), and
that she should therefore still consider the probability of Heads to
be 1/2 at this time.

Others (eg, Elga 2000) argue instead that Beauty should assign
probability 1/3 to the coin landing Heads upon wakening before
Wednesday.  One argument for this is that if the experiment were
repeated many times, 1/3 of the wakenings before Wednesday would occur
when the coin lands Heads (since Beauty wakens twice when it lands
Tails, and only once when it lands Heads).
This view can be reinforced by supposing that on each awakening Beauty
is offered a bet in which she wins 2 dollars if the coin lands Tails
and loses 3 dollars if it lands Heads.  (We suppose that Beauty knows
such a bet will always be offered.)  Beauty would not accept this bet
if she assigns probability 1/2 to Heads.  If she assigns a probability
of 1/3 to Heads, however, her expected gain is
$2\times(2/3)-3\times(1/3) = 1/3$, so she will accept, and if the
experiment is repeated many times, she will come out ahead.
Furthermore, she can work all this out on Sunday, at which time she
will wish herself to accept these bets later on.  Accepting the
argument that she should assign probability 1/2 to Heads when she is woken
later therefore requires that we accept that Beauty should override
her previous decision, even though she has no new knowledge that would
justify such a change.  This seems at least as strange as accepting
that Beauty should alter her probability of Heads from 1/2 to 1/3 even
though nothing unexpected has apparently happened.

A problem with Lewis's argument for the probability of Heads being 1/2
arises if we change the experiment so that some time after being woken
on Monday, Beauty is told that it is Monday.  If she assigned
probability 1/2 to Heads just before being told this, standard
Bayesian updating of probabilities would lead her to assign
probability 2/3 to Heads after being told it is Monday --- the
probability of Monday given Heads is~1, whereas the probability of
Monday given Tails is 1/2; this factor of two difference shifts the
previous equal probabilities for Heads and Tails so that Heads has
twice the probability of Tails.  This seems ridiculous, however, given
that at this point, Beauty knows nothing of relevance that would
distinguish her from any other person who has gotten a good night's
sleep, and is then asked to predict the toss of a fair coin.  That she
knows this coin toss will determine whether her memory is erased
\textit{in the future} should be of no relevance.  Elga's (2000)
argument for the 1/3 view is essentially to work backwards from the
assumption that Beauty should assign probability 1/2 to Heads if she
is in this situation.\footnote{Bostrom (2006) argues that Beauty
should assign probability 1/2 to Heads both before and after being
told on Monday that it is Monday.  He justifies this non-reaction to
new information on the basis that there is a shift in reference class
on being told that it is Monday.  This argument involves the use of
narrow reference classes which I argue below are untenable.  Bostrom
also argues that the betting argument can be defused by further
consideration of reference classes and indexical information, but his
reasoning applies only if repetitions of the experiment are done with
Beauty's memory being erased between each repetition.  This is
unconvincing, since one can equally well suppose that Beauty remembers
how many times the experiment was previously done.}

\subsection{\hspace*{-7pt}Relating Sleeping Beauty to SSA, SIA, and 
   FNC}\label{sec-beauty-fnc}\vspace*{-5pt}

These differing answers to the Sleeping Beauty problem can be viewed
as consequences of applying SSA+SIA or SSA$-$SIA, with the reference
class being instances of Beauty upon wakening before Wednesday.  With
SSA$-$SIA, Heads and Tails are equally likely \textit{a priori},
and Beauty's observations on wakening before Wednesday are equally
likely given Heads or Tails --- if Heads, there is only one possible
wakening, if Tails there are two, and SSA tells us that they are
equally likely, but in any case they do not differ in any identifiable
way.  Heads and Tails are therefore still equally likely once Beauty
has woken.  If we accept SSA+SIA, however, Heads have probability 1/3
\textit{a priori} (once we are in the position of Beauty after
wakening, and therefore part of the reference class), since Heads
leads to only half as many members of the reference class as Tails.
The observations on wakening are still equally likely given Heads or
Tails, so Beauty's probability of Heads after wakening should remain
1/3.

As expected from the discussion in Section~\ref{sec-fnc}, the same
conclusion as SSA+SIA is reached by applying FNC --- simply
conditioning on the full data available to Beauty upon wakening.  In
this regard, note that the even though the experiences of Beauty upon
wakening on Monday and upon wakening on Tuesday (if she is woken then)
are identical in all ``relevant'' respects, they will not be
subjectively indistinguishable.  On Monday, a fly on the wall may
crawl upwards; on Tuesday, it may crawl downwards.  Beauty's
physiological state (heart rate, blood glucose level, etc.) will not
be identical, and will affect her thoughts at least slightly.
Treating these and other differences as random, the probability of
Beauty having \textit{at some time} the exact memories and experiences
she has after being woken this time is twice as great if the coin
lands Tails than if the coin lands Heads, since with Tails there are
two chances for these experiences to occur rather than only one.  This
computation assumes that the chance on any given day of Beauty
experiencing \textit{exactly} what she finds herself experiencing is
extremely small, as will be the case in any realistic version of the
experiment.

We can see for Sleeping Beauty how introduction of the ``essences''
discussed in Section~\ref{sec-fnc} together with ``Full Indexical
Conditioning'' changes the result from that of FNC to that of
SSA$-$SIA.  We suppose that instances of Beauty on different days have
different essences.  The probability that an instance of Beauty woken
before Wednesday will have both her current memories \textit{and} her
current unique essence is the same whether she is woken once or twice
--- a second awakening doesn't provide a second chance because Beauty
on the other awakening will have the wrong essence.  Denying SIA is
equivalent to assuming that the existence of two awakenings would not
increase the probability that an instance of Beauty with her current
essence exists.  In the context of this problem, these assumptions
regarding essences seems rather contrived, but perhaps advocates of
SSA$-$SIA might maintain that essences with these characteristics are
more plausible for the other problems discussed in this paper.

\subsection{\hspace*{-7pt}Beauty and the Prince}\vspace*{-5pt}

Since I maintain that the right results are obtained by using FNC (or
SSA+SIA), I would like to provide further evidence that 1/3 is indeed
the correct probability of Heads.  For this purpose, imagine a
``companion'' of Beauty, the Prince, who is also put to sleep on
Sunday, and woken on Monday.  However, unlike Beauty, the Prince is
woken on Tuesday regardless of how the coin lands, after being
administered a drug that causes him to forget his Monday awakening.
Like Beauty, he is always woken on Wednesday and told the experiment
is over.  Beauty and the Prince will be in the same room at all times,
and will be free to discuss their situation (if both are awake).  All
this is known to both Beauty and the Prince.

If the Prince is woken before Wednesday and finds that Beauty has also
been woken, what probability should he assign to the coin landing
Heads?  Quite clearly, he should assign probability 1/3 to Heads,
since given Heads, the probability of Beauty being woken with him is
1/2, whereas given Tails this probability is 1.  Since the coin is
fair, this factor of two larger probability of what is observed given
Tails should produce a factor of two larger probability of Tails given
the observation that Beauty has been woken too. 

The Prince will tell Beauty of his conclusion.  Should Beauty
disagree, and maintain that the probability of Heads is actually 1/2?
Beauty and the Prince have the same information (and even if they
didn't, they would after discussing the situation).  Beauty will agree
that the Prince's reasoning is correct, for him.  There seem to be no
grounds for her to decide that, for her, the probability should be
different.  Furthermore, if we wish, we can disallow discussions
between Beauty and the Prince --- Beauty is intelligent enough to know
what the Prince's conclusion will be without him having to tell her.
Indeed, we can assume that the Prince is hidden from Beauty by a
curtain, as long as she knows he is there.  Does it really matter if
we go one step further and eliminate the Prince altogether?

This contradiction between the conclusions reached by Beauty and the
Prince when they both apply SSA$-$SIA is not surprising in light of
the discussion in Section~\ref{sec-companions}, since the reference
class used by Beauty does not include instances of the Prince, and
vice versa.  If instead they both use the reference class of wakenings
before Wednesday of either Beauty or the Prince, the result of
applying SSA$-$SIA changes.  Beauty will then reason on wakening
before Wednesday that if the coin landed Heads, the probability that
she will be an instance of Beauty (rather than the Prince) is 1/3,
whereas if the coin landed Tails, this probability will be 2/4, and as
a result assign probability $(1/3)\,/\/(1/3+2/4)\ =\ 2/5$ to Heads.
The Prince reasons on wakening that the probability of his being an
instance of the Prince is 2/3 if the coin landed Heads and 2/4 if it
landed Tails, so the probability of Heads is $(2/3)\,/\,(2/3+2/4)\ =\
4/7$, equivalent to odds of 4/3 in favour of Heads.  When he sees that
Beauty is also awake, his odds shift by a factor of two in favour of
Tails, producing final odds of 2/3 for Heads, corresponding to the
probability of Heads being 2/5, which matches the conclusions of
Beauty.  While these conclusions are consistent, they appear doubtful
because of their novelty, and their sensitivity to the number of
companions.

If the reference class used by Beauty and the Prince is expanded to
include all wakenings by all humans, which seems natural, the
conclusions of SSA$-$SIA change again.  If there are a large number,
$N$, of wakenings by other people, and by Beauty or the Prince on
other days, the probability that Beauty will assign to Heads upon
wakening before Wednesday will be
$(1/(N\!+\!3))\,/\,(1/(N\!+\!3)+2/(N\!+\!4)) \approx 1/3$.  The
probability the Prince assigns to Heads upon wakening will be
$(2/(N\!+\!3))\,/\,(2/(N\!+\!3)+2/(N\!+\!4)) \approx 1/2$, which will
change to $1/3$ (a shift of odds by a factor of two) when he sees that
Beauty is also awake.  These conclusions of Beauty and the Prince are
consistent, and match the conclusions obtained using FNC or SSA+SIA.

We therefore see that to obtain the answer 1/2, or any answer other
than 1/3, SSA$-$SIA must be applied with a narrow reference class.
Furthermore, the most natural narrow reference class --- instances of
Beauty alone --- cannot be used if consistency with the conclusions of
a companion is required.  The arguments based on SSA$-$SIA for the
answer 1/2 appear to not be viable.  Any such arguments in favour of
some answer other than 1/2 or 1/3 seem arbitrary, and also
unmotivated, to the extent that the original intuition in favour of
1/2 --- that Beauty learns nothing when she wakens before Wednesday
--- is violated by any other answer.

\subsection{\hspace*{-7pt}The Sailor's Child problem}\vspace*{-5pt}

The issue underlying the Sleeping Beauty problem can be further clarified
by looking at what I will call the Sailor's Child problem.  This
problem does not involve memory loss.  Without any such fantastic
aspect, we can surely hope to obtain a clear answer.

A Sailor sails regularly between two ports, in each of which he stays
with a woman, both of whom wish to have a child by him.  He is
reluctant, but eventually decides that he will have one or two
children, with the number decided by a coin toss --- one if Heads, two
if Tails.  Furthermore, he decides that if the coin lands Heads, he
will have a child with the woman who lives in the city listed first in
\textit{The Sailor's Guide to Ports}.  (He considers this fair, since
although he owns a copy of this book, he hasn't previously read it, and
so has no prior knowledge of which city comes first.)  Now, suppose
that you are this Sailor's child, and that neither you nor your mother
know whether he had a child with the other woman.  You also do not
have a copy of \textit{The Sailor's Guide to Ports}.  You do, however,
know that he decided these matters as described above.  What should
you consider to be the probability that you are his only child (ie,
that the coin he tossed landed Heads)?

The answer seems clear.  Given your ignorance regarding \textit{The
Sailor's Guide to Ports}, you should believe that if the coin landed
Heads, your mother would have been selected to have a child with
probability 1/2, whereas if the coin landed Tails, this probability
would have been~1.  This 2-to-1 ratio of probabilities for what is
observed (that you were born) given Tails versus Heads leads to the
probability of Heads being 1/3.  The probability of your having a
half-sibling is therefore 2/3.  If you have any doubts about this, due
to some idea that you should consider that you might have been the
other child, ask the opinion of your mother, who plays the role of
``companion'' in this tale.  She should have no doubts about this
reasoning.

Suppose that you later obtain a copy of \textit{The Sailor's Guide to
Ports}, and find that the city you were born in is listed before the
other port city.  With this additional information, your birth becomes
certain, regardless of the result of the coin flip.  You therefore
have no information regarding the result of this flip, and should
assign probability 1/2 to Tails, and hence also to the possibility
that you have a half-sibling.  This situation is analogous to Beauty
being told that it is Monday sometime after wakening.

Do the Sailor's Child and Sleeping Beauty problems differ in any
important way?  In the Sleeping Beauty problem, the instances of
Beauty awakening on Monday and Tuesday can be visualized as
``children'' of the Beauty who existed on Sunday.  That these
``children'' are much more closely related than are real children of
the same father seems inessential, particularly since the only
information transferred from Beauty-on-Sunday to Beauty-awakened-later
is ``common knowledge'' about the setup, such as that the coin is
fair.  In this light, we can see that contrary to many treatments in
the literature, the Sleeping Beauty problem is not really about
updating of beliefs as new information is received --- a procedure
that in any case seems dubious when actual or suspected memory loss is
an issue.

Perhaps the puzzlement concerning Sleeping Beauty is partly a
consequence of an implicit assumption that if Beauty is woken on both
Monday and Tuesday her subjective experiences will be identical ---
that she will in essence be a single sentient entity existing at two
moments in time.  As mentioned above, this assumption is not justified
by the standard description of the experiment.  And indeed, it is not
explicitly used in any of the arguments that I am aware of, but
nevertheless seems to colour thinking about the problem.  For example,
Lewis (2001, p.~171) in describing the problem says, ``\ldots the
memory erasure on Monday will make sure that her total evidence at the
Tuesday awakening is exactly the same as at the Monday awakening'',
and Elga (2000, p.~145) says, ``We may even suppose that you knew at
the start of the experiment exactly what sensory experiences you would
have upon being awakened on Monday'', which in the context of the
problem would require also assuming that these experiences are
identical to those on Tuesday (if one is woken then).  Assuming
that Beauty's subjective experiences on Monday and Tuesday are
identical converts the thought experiment from one with an only mildly
fantastic element (perfect memory erasure) to one which is arguably
impossible in principle.  It is perhaps not surprising that confusion
can then ensue.  In the Sailor's Child problem, however, no one would
assume that if the Sailor has two children their lives will be
indistinguishable; indeed, it is obvious that their lives will differ
substantially (different mothers, different cities, etc.).  That we
can imagine the Sleeping Beauty experiment happening with the
experiences of Beauty on Monday and on Tuesday differing only in much
less dramatic ways does not change the correct answer; it only makes
it harder to see.

\section{\hspace*{-7pt}The Doomsday Argument}\label{sec-da}\vspace*{-10pt}

Some versions of the Doomsday Argument, such as that of Gott (1993),
depend in essence on an unsupported intuition --- if humanity will
expand into the galaxy, with hundreds of trillions of humans being
born, isn't it rather surprising that you are among the first hundred
billion humans?  I will deal only with the version due to Leslie
(1996), which can be put in formal terms if SSA is assumed, as
discussed by Bostrom (2002).

\subsection{\hspace*{-7pt}Formalizing the Doomsday Argument}\vspace*{-5pt}

Suppose you know that you are the $r$'th human to be born (knowing $r$
approximately is good enough, but would complicate the notation).  Let $N$
be the (unknown) total number of humans who will ever be born.  If you
assume SSA with the reference class being all humans, the probability
of your observation that your birth rank, $R$, is $r$, given that
the total number of humans is $n$, can be written as 
\beq
 P(R=r\,|\,N=n) & = & \left\{ \begin{array}{cl} 1/n & \mbox{if $n \ge r$}\\[3pt]
                                                0   & \mbox{if $n < r$}
                       \end{array}\right.
\eeq
This is a consequence of SSA, since you might equally have been any of
the $n$ humans to ever exist.  Suppose that $P(N=n)$ gives your prior
belief that a total of $n$ humans will ever be born, based on
information such as your judgement of the probability that a large
asteroid will collide earth, and of the probability that a species that
has evolved by natural selection to be competitive will destroy itself
once its members acquire the technological means to do so.  Applying 
Bayes' Rule, you can obtain the posterior distribution for $N$, which 
on these assumptions
reflects what your beliefs about $N$ should be after accounting for both
your prior beliefs and your observation of $r$:
\beq
  P(N=n\,|\,R=r) & = & {\textstyle P(N=n)\,P(R=r\,|\,N=n) \over \textstyle
        \sum_{n'=1}^{\infty^{\rule{0pt}{2pt}}}P(N=n')\,P(R=r\,|\,N=n')} 
  \\[10pt]
                 & = & 
  \left\{ \begin{array}{cl} {\textstyle P(N=n)/n \over \textstyle
        \sum_{n'=r}^{\infty^{\rule{0pt}{2pt}}}P(N=n')/n'}
                                & \mbox{if $n \ge r$} \\[12pt]
                            0   & \mbox{if $n < r$}
  \end{array}\right.\label{eq-doom}
\eeq
Compare this with the posterior distribution given the 
information that $N$ is at least $r$ --- which is implied by your
observation that \textit{you} are the $r$'th human, but does not
contain any ``indexical'' information regarding yourself:
\beq
  P(N=n\,|\,N \ge r) & = &
  \left\{ \begin{array}{cl} {\textstyle P(N=n) \over \textstyle
                             \sum_{n'=r}^{\infty^{\rule{0pt}{2pt}}}P(N=n')}
                                & \mbox{if $n \ge r$} \\[12pt]
                            0   & \mbox{if $n < r$}
  \end{array}\right.\label{eq-nodoom}
\eeq

The difference between the posterior distributions for $N$ given by
(\ref{eq-doom}) and (\ref{eq-nodoom}) can be substantial.  Suppose,
for instance, that you believe that we will either destroy ourselves
soon, or if we avoid this fate, we (or our descendants) will go on to 
colonize the galaxy.  This belief would lead to a prior for $N$ that (when 
idealized a bit) can be expressed as something like the following:
\beq
  P(N=n) & = & \left\{\begin{array}{cl}
      1/2 & \mbox{if $n=10^{11}$} \\[3pt]
      1/2 & \mbox{if $n=10^{14}$} \\[3pt]
       0  & \mbox{otherwise}
  \end{array}\right.
\eeq
If you observe that your birth rank is $r=6\times10^{10}$,
conditioning on $N \ge r$ as in (\ref{eq-nodoom}) produces a posterior
distribution that is the same as the prior.  In contrast, conditioning
on $R=r$ as in (\ref{eq-doom}) produces a posterior distribution in
which the probability that humanity will colonize the galaxy is
reduced from $1/2$ in the prior to $(1/2)/10^{14}\ /\ ((1/2)/10^{11}
+(1/2)/10^{14})\ =\ 0.000999001$ in the posterior (or put another way, 
the odds in favour of colonizing the galaxy change from 1 to 1/1000).
Much greater shifts in odds are possible if galactic colonization 
is assumed to be more extensive (eg, of $10^{10}$ stars, each with population 
$10^{10}$).  A large ``Doomsday effect'' occurs even if the non-doom scenario 
involves only full utilization of our solar system.

\subsection{\hspace*{-7pt}Why the Doomsday Argument must be wrong}\vspace*{-5pt}

Although Leslie (1996) considers this shift in probabilities towards
doom to be correct, and Bostrom (2002) does not consider the argument
for it to be definitely refuted, there are several reasons for
rejecting the Doomsday Argument that I regard as convincing, even
without a detailed understanding of why it is wrong.

The biggest problem with the Doomsday Argument is that its conclusion
depends critically on the choice of reference class.  Is it all
members of the species \textit{Homo sapiens}?  If so, exactly how is
this species defined?  Or should earlier extinct species of the genus
\textit{Homo} be included?  Do the other Great Apes (gorillas,
chimpanzees, and orangutans) count?  Would they count if future
experiments showed their cognitive abilities were greater than is at
present believed --- so that predictions regarding our prospects of
colonizing the galaxy depend not just on the latest research into
possible mechanisms of interstellar travel, but also on the latest
research into whether apes can learn language?  

Such changes in the reference class could easily change the
probabilities resulting from the Doomsday Argument by a factor of ten.
The probabilities also change if your belief regarding your birth rank
changes.  Suppose we discover that the population of China from 15000
to 5000 years ago was much larger than previously supposed, so that we
estimate an extra $10^{11}$ people lived there in the past.  If you
previously estimated your birth rank as $6\times10^{10}$, you would
now estimate it as $1.6\times10^{11}$.  If your previous prior on $N$
was $P(N=10^{11})=P(N=10^{14})=1/2$, this would likely now change to
$P(N=2\times10^{11})=P(N=10^{14})=1/2$.  The result would be factor of
two improvement in the odds in favour of colonizing the galaxy (from
1/1000 to 1/500).  But is it reasonable that the latest results from
Chinese archeology should affect your beliefs in this way?

Looking into the future, how much can our descendants differ from
ourselves and still count as belonging to the reference class?
Depending on the answer, the Doomsday Argument argument might not
reduce the chances of our descendants colonizing the galaxy after all,
if we think they would do so only after they start implanting
computers in their brains, which might disqualify them from membership
in the chosen reference class.  It seems ridiculous that such an
ethical judgement of what counts as ``human'' would affect our beliefs
concerning the factual matter of whether or not our descendants will
colonize the galaxy.

The Doomsday Argument falls apart completely if one sees no reason why
SSA should be applied to individual intelligent observers (or to
briefer ``observer moments'').  Considering that the Doomsday Argument
has been discussed quite extensively, isn't the relevant unit not an
individual, but rather an intellectual community?  Due to modern
communications technology, there is currently only one intellectual
community in the world, but in the past there were many.  There will
again be many communities if civilization collapses (without humans
going extinct).  A civilization capable of colonizing the galaxy would
probably maintain communications, however.  The Doomsday Argument then
leads to the conclusion that we are more likely to colonize the galaxy
than one might have otherwise supposed.  Clearly, an advocate of the
Doomsday Argument needs to exclude this reference class, but it's not
apparent how this could be justified.

Another reason to doubt the validity of the Doomsday Argument is that
the same reasoning in slightly different contexts leads to conclusions
that also seem counterintuitive.  Consider the reverse Doomsday
Argument, in which you know your death rank, but are uncertain about
your birth rank.  For instance, you may have detected an asteroid on
collision course with earth, one sufficiently large that it will
certainly kill all humans.  You then know that you will be among the
last $7\times10^9$ humans to die.  If you accept the logic of the
Doomsday Argument, this will affect your beliefs about how many humans
have ever lived.  The effect will be even greater if you know that the
point of impact will be on the opposite side of the earth from South
Georgia Island, and you happen at present to be the only person on
that remote island.  You will then have good reason to believe that
you will be the very last human to die, as the shock wave will reach
you last.  According to the Doomsday logic, you should then begin to
seriously entertain various strange notions, such as that the
statistics you'd previously believed regarding world population were
actually faked, and that accepted historical and archeological
accounts of the past are incorrect.  

Moreover, there is no reason to consider only temporal ranks.  Suppose
you live in a small village of about 100 people high in the Himalayas.
You meet occasional visitors from elsewhere, from whom you gather that
some people live at lower altitudes.  They also inform you that your
village is the highest permanent habitation in the world.  Due to
language difficulties, however, you have obtained no clear idea of the
total world population, and think on the basis of these reports that
it might equally well be a few million or a few billion.  If you
accept the logic of the Doomsday argument, you should then downgrade
the odds of the world population being a few billion by a factor of a
thousand, compared to the possibility that the population is a few
million, on the basis that a world population of a few billion would
produce a probability (based on SSA) of your ``altitude rank'' being
around 100 (as it is) that is a thousand times smaller than that
produced if the world population is a few million.

\subsection{\hspace*{-7pt}More general Doomsday-like 
   arguments}\label{sec-gendoom}\vspace*{-5pt}

The only real role of birth or other rank in the Doomsday Argument is
to provide a definition of a set $S$ of observers in the reference
class whose size is known (at least approximately), and which you know
you are a member of.  If you know that your birth rank is $r$, then
you know that you are a member of the set of of all observers with
birth rank no larger than $r$, a set whose size you know to be $r$.
When you condition on your membership in a set $S$ whose size you known
to be $r$, and apply SSA, your beliefs about the number, $N$, of 
observers in the reference class are modified from your prior, given
by $P(N=n)$.
Specifically, since by SSA, $P(\mbox{you are in $S$}\,|\,N=n)\ =\ r/n$,
you should conclude that
\beq
  P(N=n\,|\,\mbox{you are in $S$}) & = & 
    {\textstyle P(N=n)\,P(\mbox{you are in $S$}\,|\,N=n) \over \textstyle
      \sum_{n'=1}^{\infty^{\rule{0pt}{2pt}}}
         P(N=n')\,P(\,\mbox{you are in $S$}\,|\,N=n')} \\[10pt]
                 & = & 
  \left\{ \begin{array}{cl} {\textstyle P(N=n)/n \over \textstyle
        \sum_{n'=r}^{\infty^{\rule{0pt}{2pt}}}P(N=n')/n'}
                                & \mbox{if $n \ge r$} \\[12pt]
                            0   & \mbox{if $n < r$}
  \end{array}\right.
\eeq
This parallels the Doomsday Argument of equation~(\ref{eq-doom}).

This generalization threatens to produce further counter-intuitive
results.  Perhaps advocates of the Doomsday Argument could find some
rationale for disallowing sets $S$ that are defined with too close a
reference to you (though doing so without also excluding the set of
people with birth rank no greater than yours might be a challenge).
If so, you could not conclude on this basis that humanity will soon go
extinct (and may never been very numerous) because your native
language, of which you are the last living speaker, is, has been, and
will be the native language of only a small number of people.  Many
Doomsday-like conclusions using sets that lack such an obvious
connection to the person making the inference still seem possible,
however.  Suppose, for instance, that I have little idea of my birth
rank (lacking any knowledge of archeology), but that I do know that I
am among the roughly $10^9$ humans born between the invention of
nuclear weapons and the first visit by humans to the moon.  I can
apply Doomsday-like logic (with reference class of all humans) to draw
a pessimistic conclusion regarding the total number, $N$, of humans
who ever exist.  The argument doesn't tell me whether the humans who
do exist were born before or after I was, but this just increases
pessimism with regard to the future.

As another example, suppose that you were convinced that only our
solar system contains planets and that some sort of catastrophe is
bound to wipe out humanity fairly soon.  In particular, you think
there will be no more than $10^{12}$ humans.  However, you consider it
plausible that intelligent beings may exist on Jupiter.  Moreover, you
believe that if such beings do exist, they are very numerous compared
to humans --- Jupiter is much larger than earth, supporting a very
large population at any given time, and reducing the chances that a
catastrophe could wipe out the whole population.  Suppose you think
that $10^{16}$ Jupiter beings will exist if any such beings exist.  If
you now apply SSA with the reference class of all intelligent
observers, and consider the set $S$ of human observers, you will
multiply your prior odds in favour of the existence of Jupiter-beings
by the factor $10^{12}/10^{16} = 1/10000$, since if Jupiter-beings
exist, the probability that you would be human rather a Jupiter-being
is only $1/10000$.  Many of the reasons for disbelieving the ordinary
Doomsday Argument apply to arguments of this sort as well.  For
example, how intelligent do the Jupiter-beings have to be to count as
members of the reference class?

To me, these arguments seem just as ``presumptuous'' as the
Presumptuous Philosopher's argument discussed below in
Section~\ref{sec-pp}.  Moreover, being based on SSA$-$SIA, these
Doomsday-like arguments are sensitive to choice of reference class,
unlike arguments based on FNC or SSA+SIA (as discussed in
Sections~\ref{sec-idx} and~\ref{sec-fnc}).  However, arguments of this
type are essential if SSA$-$SIA is to avoid conflict with conclusions
by companions, as discussed generally in Section~\ref{sec-companions},
and in connection with the Presumptuous Philosopher problem in
Section~\ref{sec-density}.

\subsection{\hspace*{-7pt}The Doomsday Argument with non-human intelligent
  species}\label{sec-da-et}\vspace*{-5pt}

Inclusion of intelligent species other than humans in the reference
class changes the focus of the the Doomsday argument from specifically
human hazards to hazards affecting intelligent species in general.
This is discussed, for example, by Knobe, Olum, and Vilenkin (2006).

Suppose you are sure that many intelligent species exist, have
existed, or will exist.  Direct application of the Doomsday Argument
with this reference class is then not possible, assuming you have
little idea of your birth rank among all intelligent observers.  If
you know the distribution of the total number of individuals in a
species who ever exist, there will be no Doomsday effect ---
application of SSA should lead you to consider it more likely
\textit{a~priori} that your species is one of the more numerous ones,
which cancels the Doomsday effect from knowing that your birth rank
within your species is low.  However, a Doomsday effect remains if you
are uncertain about the typical lifetime of intelligent species.  For
instance, you might be uncertain whether species that develop advanced
technology are likely to use it to destroy themselves (either by
internal conflict, or by destruction of their environment), in which
case most intelligent species will be short lived, or whether instead
such advanced technology will typically enhance a species' prospects,
in which case most intelligent species will be long lived, with
numerous individuals.  Application of SSA together with knowledge that
your birth rank within your species is low then produces a shift of
probability toward the hypothesis that most intelligent species are
short lived (most likely including yours, absent any evidence that it
is an exception).

\subsection{\hspace*{-7pt}Defusing the Doomsday Argument with 
                          SIA or FNC}\vspace*{-5pt}

Assuming SIA as well as SSA defuses the Doomsday Argument.  Your
original prior for $N$, given by probabilities $P(N=n)$, is adjusted
by SIA so that the prior probability that $N=n$ becomes
$n\,P(N=n)\ /\ \sum_{n'=1}^{\infty} n'\,P(N=n')$.  Substituting this
expression for occurrences of the form $P(N=n)$ in
equation~(\ref{eq-doom}) gives the ``no doom'' probabilities of
equation~(\ref{eq-nodoom}).

As one would expect from the discussion in Section~\ref{sec-fnc},
the same result is obtained if one uses FNC, with the argument being
even more direct.  You update your prior for $N$, with probabilities
$P(N=n)$, by conditioning on the existence of a person with your full
set of memories, which includes knowledge of your birth rank, $r$.
(Your scientific knowledge regarding what $N$ might be will have been
incorporated into the prior, however, and needn't be conditioned on at
this point.)  However, the (presumably very small) probability that a
person exists with birth rank $r$ and with all your other memories is
independent of $N$, apart from the requirement that $N$ be at least
$r$.  The posterior distribution of $N$ is therefore just the prior
for $N$ renormalized to sum to one over the range from $r$ up --- ie,
the same as in equation~(\ref{eq-nodoom}).

Similar arguments refute the form of the Doomsday Argument where there
are many intelligent species.  Assuming SIA makes it more likely that
intelligent species are usually long-lived, cancelling the doomsday
effect.  If you apply FNC, the probability that someone exists with
all your memories, including your knowledge that your birth rank
within your species is $r$, depends only on how many intelligent
species have at least $r$ individuals, not on how long-lived these
species are beyond that.

The apparent simplicity of these refutations of the Doomsday
Argument is deceptive, however.  These refutations will have no force
if assuming FNC (or SIA, if you prefer) leads to insurmountable problems
in other contexts.  This issue is explored in the next two sections.

\section{\hspace*{-7pt}Freak Observers}\label{sec-fo}\vspace*{-10pt}

Bostrom (2002) argues for SSA on the grounds that without it drawing
conclusions from empirical evidence is impossible, due to the
existence of ``freak observers''.  The problem is that in a large
universe, \textit{someone} will have made every possible observation,
regardless of what the true state of the universe is.  So knowing that
a particular observation has been made provides no evidence at all
concerning reality.  However, if you accept SSA, and know that
\textit{you} made a particular observation, you can draw useful
conclusions, as long as \textit{most} observations made by observers
in your reference class correspond (at least approximately) to
reality.  For this to work, your reference class needn't be very large
--- it makes no difference whether you use all human-like beings or
all intelligent observers, for instance --- but it must be at least a
bit larger than the set of observers with exactly your memories, since
that narrow reference class might consist only of observers who have
made the same misleading observation.

In arguing for the existence of freak observers, Bostrom mentions
possibilities such as brains in all possible states being emitted from
black holes as Hawking radiation, or condensing by chance from clouds
of gas.  However, it is perhaps too easy to say that in an infinite
universe such events must happen.  The size of the universe that is
needed to make such events likely is quite unimaginably huge, even in
comparison with the vastness of the observable universe that we have
by now become accustomed to.  It is much easier to imagine misleading
observations arising by more normal mechanisms.  Equipment failures,
unusual amounts of noise, incompetence, and fraud are all possible
reasons why an apparently definitive scientific observation might
actually be wrong.  The probability of an observation being wrong will
almost always be at least one in a billion.  Since there are quite
likely more than a billion planets in the observable universe that are
inhabited by intelligent beings (this is much less than one per
galaxy), it is likely that numerous highly misleading observations
have been made.  

I will argue that you can draw conclusions from observation despite
the existence of such misleading observations without needing SSA ---
use of FNC being sufficient --- as long as the universe is not so
large that you would expect there to be other observers with exactly
the same memories as you.  According to FNC, you should condition not
just on what you know to be the result of a scientific observation,
but also on all your other memories.  The probability that an observer
exists with all your memories, including your memory of the
observation, will be much greater if the observation is correct than
if it is incorrect, unless the universe is so large that it contains
many observers whose memories match yours in all respects other than
the result of this observation.  You will therefore be justified in
concluding that your observation likely corresponds to reality.

How big would the universe have to be for the assumption that there
are no other observers whose memories match yours to be false?  Though
I have not performed any detailed calculations, it seems to me that
the most likely way for another observer with your memories to arise
is by ordinary biological processes on a planet somewhere --- not by
bizarre mechanisms such as Hawking radiation.  To get a lower limit on
the required size of the universe, let us suppose that life-bearing
planets are common, and that biology on them is always much like it is
on earth.  Producing a duplicate of you would then require that a
species evolve that is nearly identical to \textit{Homo sapiens}, and
that an individual in this species then acquire the same memories as
you.  The human genome contains about $3\times10^9$ base pairs, for
each of which there are four possible nucleotides.  Even supposing
that only 1\% of these base pairs might differ in a similar species
(others being functionally constrained), and that only 1\% of these
differences would have a noticeable effect, we are left with
$4^{3\times10^5} \approx 10^{180000}$ possible and distinguishable
human-like individuals.  Humans have approximately $10^{11}$ neurons.
Even supposing that each neuron and its connections (which typically
number in the thousands) encodes only one bit of useful information,
there will be $2^{10^{11}} \approx 10^{30000000000}$ possible sets of
memories.  This number dominates the number of possible genomes.  For
comparison, there are roughly $10^{11}$ galaxies in the observable
universe, each containing roughly $10^{11}$ stars.  Even if all these
stars have life-bearing planets, the number of such planets in the
observable universe is only $10^{22}$.  If each contained $10^{10}$
observers, replaced over $10^{10}$ generations, the total number of observers 
would be $10^{42}$.  The universe would therefore need to be a factor of around
$10^{30000000000}/10^{42} = 10^{29999999958}$ times larger than
the portion of it that we can observe in order for there to be a good
chance that another observer exists with the same memories as you.

Large as it is, $10^{29999999958}$ is of course as nothing compared to
infinity, which some may believe describes the actual extent of the
universe (eg, Knobe, Olum, and Vilenkin 2006).  However, before
insisting that the possibility of a universe this large or infinite
should constrain our reasoning processes, we should ask what else
would change if we took this possibility seriously.  Common notions of
decisions and ethics would seem to be in serious trouble --- if
everything is bound to happen someplace, why strive for a good outcome
here?  As I argued in Section~\ref{sec-fant}, making fantastic
assumptions of this sort carries the danger that subsequent reasoning
may utilize premises that are in fact incompatible with the
assumption.

Accordingly, it seems safer to at least initially consider problems of
anthropic reasoning on the assumption of a finite (and not
ridiculously large) universe.  However, I do not dismiss the
possibility that the universe is truly infinite in spatial extent, or
that an infinity of parallel universes exist, due perhaps to the Many
Worlds interpretation of quantum mechanics being correct.  But if so,
and if the consequences are as they might superficially seem, then the
arguments concerning the puzzles of anthropic reasoning that I address
here need to go beyond what has appeared in the literature so far.
Great care would need to be taken to ensure that assumptions based on
ideas such as uniqueness of individuals do not enter in subtle ways.
I will consider infinite universes in Section~\ref{sec-cosmo}, but for
the moment, I confine myself to arguments that are compatible with a
more common sense view of the universe.

\section{\hspace*{-7pt}Presumptuous Philosophers}\label{sec-pp}\vspace*{-10pt}

The Presumptuous Philosopher problem has been interpreted (eg, by
Bostrom (2002)) as showing that SIA should not be accepted, because it
leads to an unreasonable preference for cosmological theories that
imply that the number of observers in the universe is large.  If
theory $A$ implies that the number of observers (in the chosen
reference class) is a trillion times larger than that implied by
theory $B$, SIA says that you should shift your relative prior beliefs
in theories $A$ and $B$ by a factor of a trillion.  If you judge the
two theories equally likely on other grounds (and there are no other
plausible theories), application of SIA should leave you virtually
certain that theory $A$ is true --- so certain that you would
rationally ignore almost any future evidence to the contrary.  Such
dogmatism seems intuitively unacceptable.

FNC also seems vulnerable to the Presumptuous Philosopher problem.  The
probability of an observer with exactly your memories existing
somewhere in the universe will be greater if the number of observers
who ``might have'' your memories is larger.  If your memories are
detailed enough to make the probability of their occurring small even
in the largest universe considered, the probability that an observer
with your memories exists will be directly proportional to the number
of observers, producing the same shift as for SIA. 

As an extreme, a theory that says the universe is infinite would
appear to be infinitely favoured by SIA.  FNC would also favour such a
theory over one in which the universe is only a few billion light
years in extent, though with FNC the preference for an infinite
universe would only be very large, not infinite (since once the
universe is large enough that it is nearly certain that at least one
observer with your memories will exist somewhere, sometime, further
increases in the size of the universe are neither favoured or
disfavoured).  Here, however, I will discuss the Presumptuous
Philosopher problem under the assumption that the universe is finite,
and though it may be very large, it is not so large that exact
duplicates of observers are likely to occur.  I consider first the
consequences of SSA$-$SIA versus those of SSA+SIA for theories that
differ in the density of observers, but not in the size of the
universe.  I then consider what FNC says in this situation.  Finally,
I discuss the more difficult problem of assessing theories that
predict universes of different sizes.

\subsection{\hspace*{-7pt}Theories differing in the density of observers 
                  --- SSA$-$SIA vs.\ SSA+SIA}\label{sec-density}\vspace*{-5pt}

Implicit in the simple statement of the Presumptuous Philosopher
problem is the assumption that the details of why theory $A$ predicts
many more observers than theory $B$ do not matter.  In particular,
Bostrom and \'{C}irkovi\'{c} (2003) present two versions of the
Presumptuous Philosopher problem, one in which theories $A$ and $B$
differ with respect to the size of the universe, the other in which
they differ with respect to the density of observers, and consider
(without discussing the matter) that the same conclusion (that SIA
gives wrong results) should be reached in both situations.  On the
contrary, I will argue here that if theories $A$ and $B$ agree on the
size of the universe (which both also predict is homogeneous), but
differ in that theory $A$ predicts a higher density of observers, then
the preference given by SIA to theory $A$ may be well justified. 

Throughout this section and the next, I will assume that intelligent
observers in different star systems cannot affect each other --- not
through travel, communication, or detection, nor though any
non-deliberate means --- since any contact with observers in other star
systems would provide direct evidence of the density of observers,
invalidating the discussion.  It is important to keep in mind that
this assumption is false.  Section~\ref{sec-int} discusses what we may
conclude about the density of observers given what we actually know.
The discussion below is meant to clarify the philosophical issues
involved in a simple empirical context, which does not correspond to
our actual state of knowledge.

Leslie (1996) and Olum (2002) discuss a historical instance of
theories differing with regard to the density of observers, in which
Olum argues in favour of the predictions of SIA, whereas Leslie, and
also Bostrom and \'{C}irkovi\'{c} (2003), argue against SIA.
Marochnik (1983) advanced a theory that earth-like planets can form
only around stars whose distance from the galactic centre leads them
to revolve around the galactic centre at a speed that nearly maintains
their position relative to the density waves defining the galaxy's
spiral arms.  This theory limits the possible locations of earth-like
planets to a small region, occupying a fraction $f$ of the galaxy,
with $f$ perhaps being in the range 0.01 to 0.1.  Given various other
assumptions (eg, that life formed in this region does not colonize the
rest of the galaxy), Marochnik's theory would imply that there are a
fraction $f$ fewer intelligent observers in the universe than would be
expected under an alternative ``planets everywhere'' theory, in which
there is no such restriction on where an earth-like planet can form.

As Olum notes, Marochnik's theory is therefore disfavoured by SIA ---
ie, if you accept SIA, you should reduce the probability you assign to
this theory below the probability you would have assigned to it based
on ordinary considerations.  In particular, if you thought the two
theories equally likely excluding consideration of SIA, and no other
theories are plausible, you would consider Marochnik's theory to have
probability $f/(1\!+\!f)$ after applying SIA (equivalently, the odds
in favour of Marochnik's theory shift from $1$ to $f$).  In contrast,
to someone who accepts SSA but not SIA, and who does not know the
distance of any earth-like planet from the galactic centre,
Marochnik's theory is neither favoured nor disfavoured compared to the
``planets everywhere'' theory, provided that the regions where life is
possible according to Marochnik's theory are large enough that
intelligent observers are nearly certain to have arisen in some such
region at least once in the history of the universe.  ($f$ would need
to be very small, less than about $10^{-20}$, for this condition to be
false.)

Suppose now that we are able to measure the distance of the sun from
the galactic centre, and we find that this distance is such that it
leads to the sun nearly maintaining its position with respect to the
galaxy's spiral arms, as predicted by Marochnik's theory.  In fact,
according to Marochnik's (1983) paper (though perhaps not more recent
measurements), our sun does appear to be at the required distance.  We
could easily imagine this was not known until later, however.  Should
such an observation be taken as a confirmation of Marochnik's theory?

According to all views of the matter, this observation does indeed
increase the probability that Marochnik's theory is correct.  In
particular, the odds in favour of Marochnik's theory are multiplied by
the ratio of the probabilities of this observation under Marochnik's
theory and the ``planets everywhere'' theory, which is $1/f$ (since
the probability of observing this is $1$ for Marochnik's theory and
$f$ for the ``planets everywhere'' theory).  However, if SIA is
accepted, this increase in the probability of Marochnik's theory
merely cancels the previous lowering of the theory's probability due
to its prediction that there are relatively few intelligent observers.
(Equivalently, the odds in favour of Marochnik's theory shift from $f$
to $f\times(1/f)=1$.)  The result is that Marochnik's theory, after
such a observation is made, has the same probability as it would have
had without any such observation, and without adjusting its prior
probability by applying SIA.  In contrast, someone who accepts SSA but
not SIA, and who on the basis of prior information regards Marochnik's
theory and the`` planets everywhere'' theory as equally likely, would
take the observation that our sun is in the small region of the galaxy
where Marochnik's theory predicts stars with earth-like planets are
possible as reason to increase the probability of Marochnik's theory
to $1/(1\!+\!f)$ (ie, the odds in favour of Marochnik's theory would
be $1/f$).  This preference comes about because under SSA (with a
reference class of all intelligent observers) it is unlikely that we
would be in this special place in the galaxy if intelligent observers
are found throughout the galaxy, but it is certain that we will be in
this special place if it is the only place where intelligent observers
can exist.

The device of imagining ``companion'' observers can be used to shed
light on which of these views is correct.  Suppose that in addition to
any intelligent observers who may originate on planets, intelligent
observers taking the form of complex patterns of plasma and magnetic
fields exist in the atmospheres of all stars.  Initially, let us
assume that each star harbours only around a dozen such star-beings
(and we know this).  Such star-beings therefore make up a negligible
fraction of the reference class of all intelligent observers, even if
a planet holding billions of intelligent beings is found around only
one star in a million.  The star-beings living in our sun's atmosphere
are quite willing to engage in astrophysical discussions, once we
realize they are there.  After these discussions they have the same
observational data as we do.  Will their conclusions about Marochnik's
theory agree with ours?  The answer depends on the principles by which
inference is done, as discussed below, and summarized in the tables on
the next page.

\begin{figure}[p]

\begin{center}\begin{tabular}{
|@{\hspace*{4pt}}l@{\hspace{2pt}}||cc|cc||c@{\hspace*{5pt}}c|c@{\hspace*{5pt}}c|
}
\hline
  & \multicolumn{4}{c||}{}
  & \multicolumn{4}{c@{}|}{}
\\[-8pt]
STAR-BEINGS MUCH        
  & \multicolumn{4}{c||}{Reference class only star-beings} 
  & \multicolumn{4}{c@{}|}{Reference class both} 
\\    
LESS NUMEROUS
  & \multicolumn{4}{c||}{or only planet-beings} 
  & \multicolumn{4}{c@{}|}{star-beings and planet-beings} 
\\[6pt]
  & \multicolumn{2}{c|}{SSA$-$SIA}
  & \multicolumn{2}{c||}{SSA+SIA}
  & \multicolumn{2}{c|}{SSA$-$SIA}
  & \multicolumn{2}{c@{}|}{SSA+SIA}
\\[4pt]
  & planet & star
  & planet & star
  & planet & star
  & planet & star
\\[-2pt]
  & beings & beings
  & beings & beings
  & beings & beings
  & beings & beings
\\
\hline &&&&&&&& \\[-12pt] \hline
Prior based on ordinary   &  1  &  1  &  1  &  1  &  1  &  1  &  1  &  1  \\
information &&&&&&&&  \\ \hline
Prior after adjustment    &     &     &  $f$&  1  &     &     &  $f$&  $f$\\
using SIA &&&&&&&& \\ \hline
Prior after adjustment    &  1  &  1  &  $f$&  1  &  1  &$1/f$&  $f$&  1  \\
using SSA &&&&&&&& \\ \hline
Posterior after existence &  1  &  $f$&  $f$&  $f$&  1  &  1  &  $f$&  $f$\\
of companions known &&&&&&&& \\ \hline
Posterior after location  &$1/f$&  1  &  1  &  1  &$1/f$&$1/f$&  1  &  1  \\
of sun in galaxy known &&&&&&&& \\ \hline
\end{tabular}

\vspace*{23pt}

\begin{tabular}{
|@{\hspace*{4pt}}l@{\hspace{2pt}}||cc|cc||c@{\hspace*{5pt}}c|c@{\hspace*{5pt}}c|
}
\hline
  & \multicolumn{4}{c||}{}
  & \multicolumn{4}{c@{}|}{}
\\[-8pt]
STAR-BEINGS MUCH        
  & \multicolumn{4}{c||}{Reference class only star-beings} 
  & \multicolumn{4}{c@{}|}{Reference class both} 
\\    
MORE NUMEROUS
  & \multicolumn{4}{c||}{or only planet-beings} 
  & \multicolumn{4}{c@{}|}{star-beings and planet-beings} 
\\[6pt]
  & \multicolumn{2}{c|}{SSA$-$SIA}
  & \multicolumn{2}{c||}{SSA+SIA}
  & \multicolumn{2}{c|}{SSA$-$SIA}
  & \multicolumn{2}{c@{}|}{SSA+SIA}
\\[4pt]
  & planet & star
  & planet & star
  & planet & star
  & planet & star
\\[-2pt]
  & beings & beings
  & beings & beings
  & beings & beings
  & beings & beings
\\
\hline &&&&&&&& \\[-12pt] \hline
Prior based on ordinary   &  1  &  1  &  1  &  1  &  1  &  1  &  1  &  1  \\
information &&&&&&&&  \\ \hline
Prior after adjustment    &     &     &  $f$&  1  &     &     &  1  &  1  \\
using SIA &&&&&&&& \\ \hline
Prior after adjustment    &  1  &  1  &  $f$&  1  &  $f$&  1  &  $f$&  1  \\
using SSA &&&&&&&& \\ \hline
Posterior after existence &  1  &  $f$&  $f$&  $f$&  $f$&  $f$&  $f$&  $f$\\
of companions known &&&&&&&& \\ \hline
Posterior after location  &$1/f$&  1  &  1  &  1  &  1  &  1  &  1  &  1  \\
of sun in galaxy known &&&&&&&& \\ \hline
\end{tabular}\end{center}

\vspace*{19pt}


Inferences by planet-beings and star-beings under various
assumptions.  The upper table is for the scenario in which star-beings
are much less numerous than planet-beings, the lower table when
star-beings are much more numerous than planet-beings.  Entries in the
tables are the odds, utilizing the information listed on the left for
that row and the rows above, in favour of Marocknik's theory versus
the ``planets everywhere'' theory, with $f$ being the fraction of stars in
the location where Marochnik's theory says earth-like planets are
possible.


\end{figure}

Consider first what conclusions will be drawn if humans (or other
planet-beings) use as their reference class for SSA or SIA only other
planet-beings (not star beings), and similarly star-beings use as
their reference reference class only other star-beings.  Since the two
theories under consideration make the same predictions regarding
star-beings, SSA and SIA will then have no effect on inferences by
star-beings, but may have an effect on inferences by planet-beings.
(These conclusions are shown on the left side of the top table.)

Suppose first that neither we nor the star-beings in our sun's
atmosphere know the position of the sun in the galaxy.  The
star-beings will then take our existence on earth as evidence against
Marochnik's theory, since according to that theory, earth-like planets
occur only in a special region of the galaxy, and (like us) the
star-beings have no reason to think that the sun is in this special
region.  If they thought Marochnik's theory and the ``planets
everywhere'' theory were equally likely before knowing that humans
exist, they will think Marochnik's theory has probability
$f/(1\!+\!f)$ after (ie, the odds in favour of Marochnik's theory will
be $f$).

Suppose now that the distance of the sun from the galactic centre is
measured, and found to be such that the sun nearly maintains its
position with respect to the spiral arms.  The star-beings will now no
longer take our existence on earth as evidence against Marochnik's
theory, since given that the sun is in this special place, it is not
surprising that it has an earth-like planet.  However, Marochnik's
theory is not made any more probable than the ``planets everywhere''
theory by this observation.  As discussed above, after finding that
humans exist, the star-beings' odds for Marochnik's theory being true
would be $f$.  The probabilities for the subsequent observation that
the sun is in the special position required by Marochnik's theory are
$1$ if Marochnik's theory is true, and $f$ if the ``planets
everywhere'' theory is true.  The odds for Marochnik's theory shift by
the ratio of these probabilities, leaving the odds at $f\times(1/f)=1$.

To see why it is reasonable that the star-beings would not favour
Marochnik's theory in this situation, consider for comparison how you
would evaluate a theory that a certain fish occurs only in acidic
lakes, versus the contrary theory that it occurs in all lakes.  The
observation that this fish is present in a nearby acidic lake tells
you nothing about which of these theories is true --- you need to look
in a non-acidic lake to obtain any relevant data.\footnote{There may
be a curious order dependence of intuitions here.  If you first find
that this lake is acidic, and then find that it contains this fish,
you will certainly reason as above.  But if you first discover that
the lake contains this fish, and only later find that the lake is
acidic, you might be tempted to take this observation as confirmation
that the fish is found only in acidic lakes, particularly if acidic
lakes are rare.  The difference appears to derive from a heuristic for
avoiding self-deception --- predicting that the lake is acidic after
finding the fish in it is psychologically risky, in that you might be
wrong.  So if you actually make such a prediction, your \textit{prior}
probability that the fish is found only in acidic lakes must really be
high.  In contrast, no such stark confrontation with reality occurs
when you know the lake is acidic before looking for the fish (which is
therefore likely to be found regardless of which theory is correct).
It might then be easier to deceive yourself regarding your true prior
beliefs.  So, taking into account your own capacity for
self-deception, your conclusion may depend on the order of
observations.  Here, however, I assume that we are not prone to
self-deception of this sort.}

The conclusions of the star-beings regarding Marochnik's theory are
the same as we would reach by applying SIA.  I will call these the
``non-anthropic'' conclusions, since they are also what one would
obtain by ignoring both SSA and SIA, as some observer outside the
universe would do (just as with the fish example above).

Before considering this reasoning by companion star-beings as
supporting the non-anthropic conclusions, however, we should consider
that the star-beings might apply SSA and/or SIA with the reference
class of all intelligent observers, including both us and them.  (The
\mbox{resulting} conclusions are shown on the right of the top table.)
If the star-beings accept SIA, they will adjust their prior to
disfavour Marochnik's theory, since it predicts many fewer intelligent
observers, thereby reducing the odds in favour of Marochnik's theory
to approximately $f$.  (Recall here that the star-beings form a
negligible fraction of all intelligent observers.)  However, if they
now apply SSA in conjunction with their observation that they are
star-beings, this effect is undone, since the more planet-based
intelligent observers there are, the smaller the odds that one is a
star-being.  (Recall that the number of star-beings is the same for
both theories).  The odds in favour of Marochnik's theory therefore
shift by the factor $1/f$, to $f\times(1/f)=1$.  Subsequent reasoning
proceeds to the non-anthropic conclusions just as above --- the
observation of humans on earth decreases the probability of
Marochnik's theory (the odds in favour decline to $f$), and the
subsequent measurement showing that the sun is in the special place
where earth-like planets are possible restores this probability to
its original value, but no higher (ie, at this point, Marochnik's
theory is neither favoured nor disfavoured).

On the other hand, if the star-beings accept SSA but not SIA, and use
the reference class of all intelligent beings, their conclusions will
match those of humans who also accept SSA but not SIA (with the same
reference class).  A doomsday-style effect occurs, of the sort
discussed in Section~\ref{sec-gendoom}, in which the star-beings
initially adjust the probability of Marochnik's theory upwards, since
this theory makes it more likely that one is a star-being, like them,
rather than a planet-being.  Observation of humans on earth causes the
probability of Marochnik's theory to go down again (since it says that
earth-like planets are rare).  Finding that the sun is in the special
place where Marochnik's theory says earth-like planets are possible
then raises the probability of Marochnik's theory back up again, to
the point where it is favoured by the factor $1/f$ over the ``planets
everywhere'' theory.

This consideration of companion star-beings therefore does not
definitively refute the methodology of accepting SSA but not SIA ---
provided humans and star-beings both use SSA$-$SIA with the reference
class of all intelligent beings, their conclusions are the same.  This
is expected from the (only slightly different) general discussion in
Section~\ref{sec-gendoom}.

However, this consideration of companions does undermine the claim
that the Presumptuous Philosopher problem makes SIA implausible (at
least on the basis of this scenario).  To review, the claim is that
although SIA defuses the counter-intuitive Doomsday Argument, it does
so at the cost of producing an equally counter-intuitive Presumptuous
Philosopher effect.  Here, however, we see that the conclusion of the
Presumptuous Philosopher in this scenario is implausible only if you
accept arguments of the Doomsday type --- such as would lead the
star-beings to favour Marochnik's theory prior to any observations
solely on the basis that it makes it more probable that an observer
will be a star-being like themselves.  An advocate of SIA who rejects
Doomsday-type arguments will therefore be untroubled by this instance
of the Presumptuous Philosopher problem.

We might alternatively suppose that, rather than each star having only
a dozen star-beings, each star instead has trillions of them, so they
are vastly more numerous than planet-based beings, under any theory.
(Conclusions on this assumption are shown in the bottom table.)  In
this situation, anthropic reasoning has no effect on the conclusions
of the star-beings, who will reach the non-anthropic conclusions
regardless of whether they accept SSA and/or SIA, using any plausible
reference class.  SIA also has no effect for human observers in this
situation, provided that they use the reference class of all
intelligent observers (not all humans, or all planet-based observers).
However, SSA has the effect, for humans, of decreasing the probability
of Marochnik's theory, since under this theory, one is less likely to
be a planet-based being than one is under the ``planets everywhere''
theory.  This decrease is undone once the sun is found to be in the
special place where planets are possible under Marochnik's theory.
The results match the non-anthropic conclusions of the star-beings,
once they know of the existence of humans, so this scenario does not
resolve any issues.  However, if the humans used the reference class
of all planet-based intelligent observers (excluding star-beings),
their conclusions using SSA$-$SIA will not match those reached by the
star-beings.

Here again, SSA$-$SIA can produce conclusions consistent with those of
companion observers provided a reference class including both is used.
Notice, however, that the conclusions found using SSA$-$SIA in this
scenario with numerous companions are different from those found when
the companions were less numerous.  Applying anthropic reasoning based
on SSA$-$SIA leads to conclusions regarding planet-based observers
that depend on how many non-planet-based observers exist, even though
this information would appear to be irrelevant.

In another scenario, we might have two theories, one of which says
that intelligent observers are common, whereas the other says that
they are randomly distributed at a low density (but high enough that
it is likely that at least one intelligent species does exist).  In
contrast with Marochnik's theory, in this scenario there is no
possibility of discovering that the sun is in a special life-bearing
region.  Using SSA$-$SIA, the prior probabilities for the two theories
depend only on normal considerations of plausibility, and our
observation that we exist does not change these probabilities.  Using
SSA+SIA, the prior probabilities are adjusted in favour of the theory
that intelligent observers are common, and these probabilities are
again unchanged by the observation that we exist.  Consideration of
star-beings as possible companion observers produces results analogous
to those discussed above for Marochnik's theory --- consistency
requires that the reference class include all intelligent observers,
and the conclusions using SSA$-$SIA, but not SSA+SIA, depend on the
number of companion observers.

In summary, consideration of companion observers provides strong
evidence that in scenarios where the density of observers varies, the
conclusions found using SSA+SIA are correct, whereas those found using
SSA$-$SIA are not.  Certainly, the conclusions of SSA$-$SIA regarding
Marochnik's theory seem quite acceptable.  Lingering unease may
remain, however, when the shift in prior odds produced by SIA is not
the factor of 10 to 100 that occurs with Marochnik's theory, but
rather a factor of a trillion or more, which one can imagine could
occur with some other theory.  I will consider such cases of extreme
``presumption'' below, when discussing FNC, which I see as the more
principled, even if largely equivalent, alternative to SSA+SIA.

\subsection{\hspace*{-7pt}Theories differing in the density of observers 
                --- result of applying FNC}\label{sec-density-fnc}\vspace*{-5pt}

To apply FNC, you multiply your prior odds for one theory over
another, based on ordinary scientific evidence, by the ratio of the
probabilities that these theories give for someone to exist with your
exact memories (excluding your scientific knowledge that contributed
to the original prior odds, but including any additional scientific
observations).  As discussed in Section~\ref{sec-fnc}, the results of
applying FNC are much the same as those of applying SSA+SIA.  Thinking
in terms of FNC avoids the need to specify any reference class of
observers, however, and clarifies the issues that are involved in this
type of reasoning.

Here is a superficial analysis of how FNC applies to
Marochnik's theory versus the ``planets everywhere'' theory, which was
discussed in Section~\ref{sec-density} using SSA$-$SIA and SSA+SIA.
If a total of $C$ earth-like planets exist according to the ``planets
everywhere'' theory, and each has some very small probability,
$\epsilon$, of producing someone with your exact memories, then the
probability of your existing according to the ``planets everywhere''
theory is $\epsilon C$ (assuming, as I am, that this is much less than
one).  According to Marochnik's theory, the number of earth-like
planets will be smaller by the factor $f$, and hence the probability
that you exist will be only $\epsilon f C$.  Assuming equal prior
probabilities for the two theories, your odds in favour of Marochnik's
theory should be $f$, if this is all you know.  Suppose you then make
a reliable observation that the sun is in the special region where,
according to Marochnik's theory, earth-like planets can form.  Your
odds in favour of Marochnik's theory should then change to~1 (ie, the
two theories should become equally likely), since the chance that
someone exists with your memories --- including your memory of
observing that the sun is in this region --- is the same for both
theories.  The extra earth-like planets that exist according to the
``planets everywhere'' theory do not increase the probability that you
will exist, since a being on a planet outside the special region is
very unlikely to remember having observed that they are in the special
region, and also have your other memories.  (See Section~\ref{sec-fo}
for consideration of the possibility of false observations.)  These
conclusions are the same as those found in Section~\ref{sec-density}
using SSA+SIA.

Although I believe that this superficial analysis produces the correct
result, it conceals a number of subtleties.  One is that, contrary to
what is implicitly assumed in the argument above, earth-like planets
in other galaxies could not really produce someone having your
memories, assuming that you have seen the numerous photos of
distant galaxies that most people have seen.  From a planet in a
different galaxy, these galaxies would be viewed from a different
angle, be much larger or smaller, or be obscured.  (As discussed
earlier, I assume for the moment that the universe is not so large
that one would expect to find another galaxy whose views of other
galaxies are just by chance virtually identical to ours.)  Because of
this, it is actually essential to the argument that Marochnik's theory
applies not just to most galaxies, but to our galaxy in particular.
If our galaxy were a rare exception in which earth-like planets were
possible everywhere even according to Marochnik's theory (and you knew
this), applying FNC would not change the probabilities of the two
theories, since the probability that you would exist in our galaxy
(where you clearly are) would be the same for both theories.  Careful
application of SSA+SIA would also lead to this conclusion, but only
because the shift in odds away from Marochnik's theory that is
produced by SIA is cancelled by applying SSA, taking account of the
greater probability of being in an exceptional galaxy if the other
galaxies are less populated.  This seems to me to be a rather perverse
way of reasoning to the correct conclusion, however.

Furthermore, it is questionable whether a planet in a distant part of
our galaxy --- likely differing from earth in elemental abundances (as
determined by local supernovae), in the local density of stars, in
cosmic ray intensities, and in the view of the Milky Way in the night
sky --- would have even a tiny chance of producing life that evolves
in just the way it has on earth, and of then producing an individual
with your memories.  As an analogy, suppose you are given detailed
photographs of a house, which you are told is either in India or in
Canada, and are asked to guess in which of these countries the house
is located.  If you are knowledgeable about architectural styles and
construction practices in the two countries, it is quite likely that
you would be able to tell which country the house is located in.  Of
course, someone with less knowledge might not be able to tell where
the house is located.  Similarly, someone with sufficient knowledge of
our galaxy would likely be able to tell where in our galaxy earth is
located, without the need for any explicit measurement of location.
However, if you lack sufficient knowledge, you will not know where in
the galaxy earth is located without an explicit measurement, and so,
\textit{as far as you know}, any earth-like planet in the galaxy might
have produced someone with your memories, and the more such earth-like
planets there are, the greater the chance that you will exist.  Put
another way, the narrow region where earth-like planets are possible
according to Marochnik's theory leads to a restricted range of
possible characteristics of these planets and their inhabitants.
Since you do not know what this range is, you do not know whether or
not the characteristics of earth and humanity are included.  The
possibility that they are not reduces the chance of your existing if
Marochnik's theory is true.  As this example makes clear, the
probabilities used in FNC may reflect your ignorance, rather than the
operation of some random physical process.

When applying FNC, it is clear that some ``presumptuous'' conclusions
that may appear to follow from SIA are not actually problematic.
Consider, for example, the theory that all bacteria are intelligent
beings.  You may regard this theory as unlikely, and assign it a low
prior probability.  However, there are approximately $10^{21}$ times
as many bacteria as humans on earth (Whitman, Coleman, and Wiebe
1998).  Similar ratios for analogous organisms presumably hold on
other earth-like planets.  According to SIA, we should therefore
increase the prior odds that bacteria are intelligent by a factor of
$10^{21}$, which may well make the theory highly probable despite its
prior implausibility.  However, if you apply FNC rather than SSA+SIA,
the probability of this theory will not be increased --- whether
bacteria are intelligent or not has no effect on the probability that
you will exist with all your memories, since you are not a bacterium.
One does in fact reach this same conclusion in the end using SSA+SIA,
since the huge increase in the theory's probability from applying SIA
is cancelled by an equally huge decrease from the low probability of
an intelligent observer being human if the theory is true.  As was the
case above, however, such reasoning based on SSA+SIA seems rather
contorted, even if the right answer is obtained, compared to the
straightforward application of FNC.  Note also that according to
SSA$-$SIA, you should decrease your (presumably already low) odds in
favour of bacteria being intelligent by a factor of $10^{21}$, on the
grounds that if they were intelligent, you would likely be a
bacterium.  In this scenario it is SSA$-$SIA, not SSA+SIA or FNC, that
could be accused of presumption.

Scenarios more troubling for an advocate of FNC can be imagined,
however.  Suppose you have calculated that the number of earth-like
planets in the galaxy is about one thousand, on the basis of what you
believe to be the correct mechanism of planet formation, and assuming
that Newton's theory of gravity is an adequate approximation.  It
occurs to you that perhaps Einstein's theory of gravity would give
different results.  You think the chances of this are only about 9\%
(odds of about 1/10), since Newton's theory is usually a good
approximation, but you decide nevertheless to redo the calculation
using Einstein's more accurate theory.  This new calculation says that
the number of earth-like planets in the galaxy is about one billion
--- a million times more than found with the Newtonian calculation.
You judge that mistakes in such calculations happen about 10\% of the
time (at least without further checking, which you haven't done yet),
so the probability of getting a divergent result such as you obtained
if Newton's theory is actually an adequate approximation is 10\%
(since a mistake would need to be made), whereas the probability of a
divergent result if Newton's theory is not adequate is about one
(since a correct result would differ, and an error would also be
fairly likely to produce a different result).  Using ordinary
reasoning, the result of this calculation should therefore lead you to
multiply by 10 the original odds in favour of the Newtonian
calculation being wrong, which produces odds of about 1.  So at this
point, you would consider that the number of earth-like planets in the
galaxy is equally likely to be one thousand or one billion.

However, if you now apply FNC (or SIA), you will increase the odds in
favour of the Newtonian calculation being wrong by a factor of a
million, since the calculation using Einstein's theory leads to a
factor of a million more earth-like planets, with a corresponding
increase in the probability of someone with your memories existing.
This extreme certainty seems presumptuous, particularly when you
haven't even checked your calculation yet.

A hint at resolving this problem comes from considering a scenario
that is similar except that the calculations are not of the number of
earth-like planets, but rather of a numerical quantity that has been
precisely measured by experiment.  If your calculation using
Einstein's theory produces a very good fit to the experimental data,
you might indeed be highly confident that it is correct, even before
checking it.  When calculating the number of earth-like planets,
however, no precise target number is matched --- FNC and SIA just say
that bigger is better, up to whatever limit is imposed by other
observations.

Accordingly, even if you accept that the number of earth-like planets
with human-like observers must be large, there is no necessity that
this number be large \textit{for this particular reason}.  There may
be many ways that the probability of your existing could be
increased other than by increasing the number of earth-like planets
--- for instance, by a higher chance of life developing on each
planet, or a higher chance that once life develops it produces an
intelligent species.  If there were no upper limit, FNC or SIA would
just favour all of these, but if there is a limit on the density of
intelligent observers, only a limited number of these factors can
strongly favour more human-like observers, reducing the probability
that any one of them in particular does.

At this point, recall the assumption stated at the beginning of
Section~\ref{sec-density} --- that intelligent observers in different
star systems have no effect on each other.  If this is true, we can
have no bound from observation on the density of intelligent
observers, and the considerations just discussed will not reduce the
excessive certainty of the conclusions from FNC and SIA.  However, if
there is a limit on the density of observers, FNC (and SIA) need not
produce unreasonably certain belief in particular theories, such as
that the hypothetical Newtonian calculation above is wrong.  We in
fact know that intelligent species may possibly interact.  The
implications of this are discussed in detail in Section~\ref{sec-int},
and show that the actual effects of FNC are not what what one would
think from a simplistic consideration of the Presumptuous Philosopher
problem.

\subsection{\hspace*{-7pt}Theories differing in the size of the 
                          universe}\label{sec-varsize}\vspace*{-5pt}

I will now consider the effect that FNC has on the probabilities of
theories that differ in the size of universe that they predict.  The
effects of applying SSA+SIA are similar, but I will omit the details
of this here.  As above, I will assume that all theories predict a
universe of a finite size, which moreover is not so large that you
would expect another observer with exactly your memories to exist.  I
will also assume that all theories predict a homogeneous universe, in
which intelligent observers arise with some density.  In practice,
different theories might well predict both different sizes of the
universe and different densities of observers, but for simplicity, I
will assume here that all theories predict the same density, so that
the total number of observers is simply proportional to the size of
the universe.

With these assumptions, a Presumptuous Philosopher effect can easily
arise.  Suppose theory $A$ says that the universe contains $10^{24}$
galaxies, whereas theory $B$ says that it contains only $10^{12}$
galaxies.  If these theories appear equally likely on ordinary
grounds, application of FNC will lead you to consider theory $A$ to be
a trillion times more likely than theory $B$, since it is a trillion
times more likely that someone with your memories will exist somewhere
if theory $A$ is true.\footnote{Here ``somewhere'' could be anywhere
in the universe.  In Section~\ref{sec-density-fnc}, where the
large-scale features of the universe could be considered fixed (since
the theories did not differ in this regard), I pointed out that
someone with your memories could exist only in our galaxy.  This is
not true in this context, where a galactic neighborhood matching what
we see is more likely to exist if the universe is large than if it is
small.}  Note that this implication of FNC holds regardless of any
details of where and how often human-like or other intelligent
observers might or might not arise --- as long as these details are
the same for both theories, the theory producing a bigger universe
will also produce a greater probability that a being with your
memories exists, in direct proportion.\footnote{Theories in which the
universe changes more slowly, and so stays longer in something
resembling the state you currently observe, might also be favoured.  I
will not elaborate on this possibility, however, but merely assume
that the universe evolves at the same rate in all theories.}

This factor of a trillion preference for theory $A$ seems unreasonable
to most people.  FNC will produce even greater degrees of certainty in
favour of theories predicting even bigger universes, up to odds of
$10^{30000000000}$ or more, before the assumption of no duplicate
observers breaks down (see Section~\ref{sec-fo}).  Unlike the
situation with theories differing in the density of observers, there
seems to be no plausible story involving companion observers that
would provide any support for this result of FNC --- if the star-being
companions of Section~\ref{sec-density} employ FNC in this situation,
they will come to the same extreme conclusion, but will also be
subject to the same intuition that this is unreasonable.

Olum (2002) offers a way of avoiding the extreme preference for larger
universes produced by FNC (or in his case, SIA) --- reduce the prior
probability of theories in proportion to the size of universe they
predict.  In the example above, if we assign theory $A$ a prior
probability a trillion times less than that assigned to theory $B$,
then after the multiplication by a factor of a trillion that comes
from applying FNC, the final odds in favour of theory $A$ will be~1
(ie, $A$ and $B$ will be considered equally likely).  This seems
rather contrived, but it does raise a crucial question --- how should
prior probabilities for cosmological theories be assigned?

For many theories, we can assign well-justified prior probabilities
based on a wealth of background knowledge.  Consider, for example,
theories regarding where eels spawn.  We can assess their prior
plausibility using our knowledge of the behaviour of other fish, as
well as our knowledge of related matters, such as ocean currents.  In
other situations, our prior beliefs will have a less detailed basis,
but will at least incorporate various common-sense constraints.  The
background knowledge we use to set priors will itself be based partly
on deeper prior beliefs.  If we could trace the origins of our beliefs
back far enough, we would presumably find some genetically-determined
prior biases, that result from natural selection.  When assessing
theories of biology, geology, macroscopic physics, or other
earth-bound phenomena, knowing that our prior beliefs have this
ultimate origin is reassuring --- it gives us some reason to think
these prior biases are well founded.

It is difficult to see, however, why natural selection should have
provided us with genetically-determined biases suitable for assigning
prior probabilities to cosmological theories.  Suppose, for example,
that the crucial difference between theories $A$ and $B$ above is that
$A$ says space is flat, with the topology of a torus, whereas $B$ says
space is positively curved, with the topology of a sphere.  (Assume
that for some theoretical reason not in dispute, the torus must have a
trillion times the volume of the sphere.)  Canceling the effect of FNC
by deciding that the spherical universe should have prior probability
a trillion times greater than the toroidal universe seems quite
arbitrary, but deciding that they should have equal prior
probabilities is really just as arbitrary.  We simply have no basis
for any prior beliefs regarding the topology of the universe.

A further difficulty is that we have no firm basis for excluding
``extraneous'' multiple universes.  Suppose the advocates of theory
$B$ above modify it to produce theory $B^*$, which is just like $B$,
except that rather than postulating the existence of only the single
universe we observe, it claims that there are a trillion similar
universes, which differ only in the actual results of random physical
processes.\footnote{Of course, this will not work if theory $B$ is
deterministic, but there will likely be ``pseudo-random'' aspects of
any theory, sensitive to slight changes in parameters or initial
conditions, that would again allow for a multiplicity of similar, but
not identical, universes.}  The probability of someone with your
memories existing in any of these universes is now the same as for
theory $A$.  This maneuver may appear unaesthetic, at least if these
trillion universes have no possibility of interacting, but compared to
the previous situation with odds of a trillion against theory $B$,
lack of elegance is a minor problem.

Is rationally assigning prior probabilities to cosmological theories
simply impossible?  Perhaps not entirely.  Sometimes, the theories
being compared all assume the same fundamental physical laws, but
represent different calculations of the consequences of these laws.
Two theories of galaxy formation, for example, are likely to assume
the same laws for gravitation and other forces, and may also assume
the same initial conditions from the big bang.  If so, the theories
can be seen a making different approximations to a single mathematical
result, whose exact computation is infeasible.  We have some
experience with mathematical approximations, and so have some basis
for assigning prior probabilities to which (if any) of these theories
is correct.

Because of the lack of clarity surrounding these issues, I see no
clear grounds for rejecting FNC or SIA on the basis of their
supposedly counterintuitive consequences regarding theories with
differently-sized universes.  Greater clarity might be obtained by
considering examples of actual cosmological theories that predict
universes of different (finite) sizes, but I am not aware of any such
examples.  Most current cosmological theories favour a universe, or
universes, of infinite size, as I discuss below in
Section~\ref{sec-inf}.

\section{\hspace*{-7pt}Anthropic arguments in
         cosmology}\label{sec-cosmo}\vspace*{-10pt}

I conclude by applying FNC to some interesting problems in cosmology,
some of which will also help further clarify the general issues
involving FNC, SSA, and SIA.

\subsection{\hspace*{-7pt}How densely do intelligent 
                          species occur?}\label{sec-int}\vspace*{-5pt}

I start with what can be seen as a continuation of the discussion
regarding the density of observers in Sections~\ref{sec-density}
and~\ref{sec-density-fnc}, except that I will now apply FNC in
conjunction with our actual empirical knowledge, which does not
correspond to the earlier assumption that there is no possibility of
intelligent species interacting.  

Our knowledge of astronomy and technology leads us to believe that if
intelligent extraterrestrials exist, it would probably not be
tremendously difficult for them, or at least their robotic probes, to
visit earth.  (We ourselves are likely to have this capability within
at most a few hundred years, unless our technological civilization
collapses.)  Interstellar travel is likely to be costly, of course,
and will certainly require patience, due to the speed of light limit.
Many extraterrestrials may decide not to undertake such exploration.
But if a large number of extraterrestrial species exist, it seems
certain that at least a fair number will explore neighboring stars.
Sometimes exploration will be followed by colonization, producing a
sphere of habitation that expands at perhaps 1\% of the speed of
light.  In this manner, a single species could reach most of the
galaxy in around 10~million years, which is a small fraction of the
galaxy's age.

Radio communication with extraterrestrial civilizations is much easier
than travel, and is well within our current capabilities (at least if
transmissions are directed at a particular star).

Despite this, we do not currently observe any extraterrestrials, nor
do we see any evidence that they have been in our vicinity in the
past.\footnote{Readers who believe they have observed
extraterrestrials may of course apply FNC themselves, and reach
different conclusions than I do here.}  This conflict with
expectations has been called ``Fermi's Paradox'', and has prompted
many attempts at explanation, summarized in a review by Brin (1983).
The paradox seems even more severe if you consider FNC (or SIA) to be
a correct principle of inference, since it seems there would then be a
further bias in favour of a high density of intelligent
extraterrestrial species (of the sort who ``might have'' produced an
observer with your memories).

The application of FNC to this problem is actually more subtle,
however.  For someone with your memories to exist, it is necessary not
only for a suitable planet to exist, but also for the subsequent
evolution of an intelligent species on that planet to proceed without
disturbance by other intelligent species.  Once someone like you
is produced, they must remain unaware of contact with any other
intelligent species.  There is therefore a tradeoff.  If earth-like
planets are common, if life arises easily on each planet, if
intelligence species are likely to evolve, and to develop a
technological civilization, the existence of someone with your
memories will be more likely, \textit{provided} there is no
interference by some other intelligent species.  But these same
factors increase the probability that such an interfering species will
exist.

A realistic analysis of this situation would be complex, as can be
seen from earlier related work on interstellar colonization, such as
that of Hanson (1998).  I will consider only a fairly simple and
abstract model intended to show the general nature of the tradeoff
described above.  This model has three components.

First, suppose that the mechanisms of galaxy formation are known, and
that the pattern of stars in our galaxy is beyond the influence of
intelligent life.  Someone with your memories, which include your
memory of the night sky, can then only have arisen on a planet of our
sun, at the current time.  Suppose, however, that we have various
theories regarding planet formation, the origin of life, the evolution
of intelligence, and the development of technological civilization.
Any particular combination of theories will produce some (tiny)
probability, $p$, that an individual with your memories will arise,
\textit{assuming} that this development is not interfered with by a
species from elsewhere.  Note that you don't know $p$, since you
don't know which theories are true, though you have prior
probabilities for them based on ordinary scientific knowledge.

Second, suppose the probability that an intelligent species with our
level of science and technology will arise in a region of size $dw$
around spacetime point $w$ is $p\, M(w)\, dw$, where $M(w)$ is a known
function giving the relative densities of intelligent species
originating at different times and places.  $M(w)$ will be zero
outside of galaxies, and at times too early for life to have
developed.  Making the probability of such an intelligent species
arising elsewhere be proportional to $p$ incorporates the assumption
that the unknown factors that influence the probability of your
existence are the same as those that influence the probability of
other intelligent species arising.

Third, suppose there is a known function, $A(w)$, and an unknown
factor, $f$, such that $f A(w)$ is the probability that a species
arising at spacetime location $w$ will destroy the possibility of
someone with your memories existing --- either by colonizing earth and
thereby preventing the development of humans, or by simply making its
existence known to you before the present time, contrary to your
actual memories.  Assuming influences are limited by the speed of
light, $A(w)$ will be zero if $w$ is outside your past light-cone.
The factor $f$ will depend on how stable technological civilizations
are, how easy interstellar travel is, and how often intelligent
species are motivated to communicate, explore, and colonize.  Assume
you have prior distributions for these factors, and hence also for $f$.

We can now write the expected number of other species that interfere with
your existence as follows:\vspace*{-10pt}
\beq
   \int\, f A(w)\, p\, M(w)\, dw & = & fp\,\int A(w)M(w)\, dw \ \ =\ \ fp V
\eeq
where $V=\int A(w)M(w)\, dw$.  Suppose that either $fpV$ is small (of order 1 
or less), or any interference from distant spacetime points is largely 
independent, so that the distribution of the number of interfering species 
will be approximately Poisson.  The probability that no species interferes
will then be $\exp(-fpV)$, and hence the probability that someone with your 
memories exists will be\vspace*{-10pt}
\beq
   P(\mbox{someone like you exists}) & = & p\, \exp(-fpV)
\label{eq-you-given-fp}
\eeq
This is maximized when $p\,=\,1\,/\,fV$, corresponding to $fpV=1$.
Thus we see that although FNC favours as large a value of $p$ as
possible when there are no interactions between species, this is not
true when interactions such as those modeled here exist, thereby
justifying the comments at the end of Section~\ref{sec-density-fnc}.

The Fermi Paradox now seems unsurprising.  If the expected number of
other intelligent species to influence earth, which is equal to $fpV$,
is around one, we should not be especially surprised that we have not
seen evidence of any other species.  We still have no specific
explanation of what factors are responsible for this, however.  In the
other direction, discovery of another intelligent species would also
not be surprising, especially if we looked somewhat more widely than
the region where a species would have influenced us without effort on
our part.

Further analysis requires some assumptions about your uncertainty
regarding $p$ and $f$.  If many unknown factors affect $p$ and $f$, in
a multiplicative fashion, it may be reasonable (due to the Central
Limit Theorem) to suppose that $\log(p)$ and $\log(f)$ have Gaussian
prior distributions.  It is also plausible that $p$ and $f$ are
independent, \textit{a~priori}.

Multiplying this prior density for $\log(p)$ and $\log(f)$ by the
probability that you exist for given values of $p$ and $f$, from
equation~(\ref{eq-you-given-fp}), and renormalizing, gives the
posterior joint probability density for $\log(p)$ and $\log(f)$.  This
density is not analytically tractable, but is easily displayed
graphically by means of a sample of points, as shown in the
accompanying figure.  Note that the numerical magnitude of $p$ depends
on exactly how detailed your memories are, and hence is probably not
of much interest.  The scale of $f$ is arbitrary, since it can be
compensated for by a change in the scale of $A(w)$, and hence of $V$.
It is convenient to set the scale of $f$ so that the mean of
$\log(p)+\log(f)$ is zero.  The parameters of interest are then the
value of $V$ and the standard deviations of $\log(p)$ and $\log(f)$.
The top-left plot shows a sample of 500 points from
the prior with standard deviation for $\log_{10}(p)$ of 1.25, giving a
95\% central interval for $p$ that spans a range of $10^5$, and
standard deviation for $\log_{10}(f)$ of 0.75, giving a 95\% central
interval for $f$ that spans a range of $10^3$.

\begin{figure}[p]

\vspace*{-25pt}

\hspace*{10pt}\includegraphics{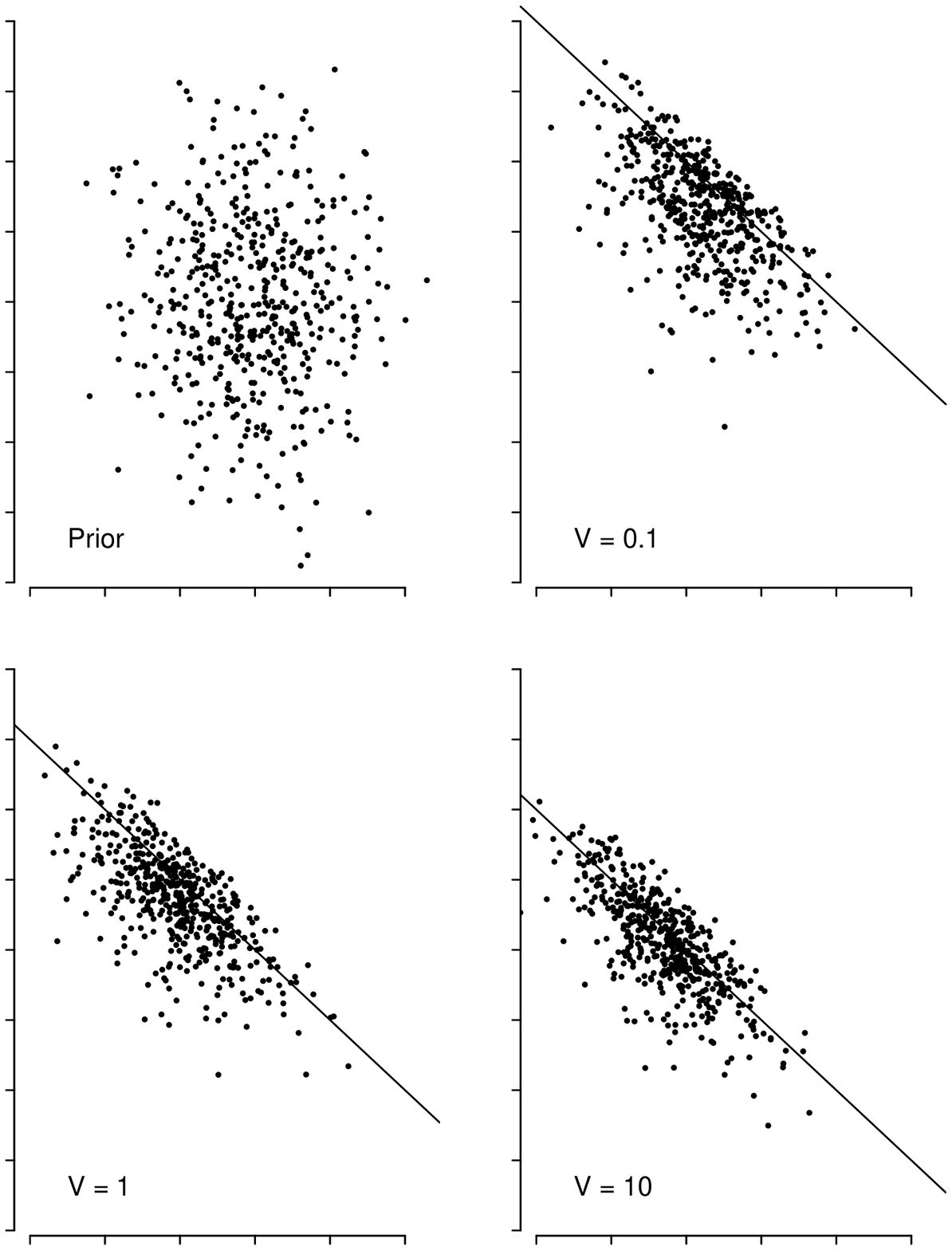}

\vspace*{12pt}

Plots of prior and posterior distributions for $\log_{10}(f)$
(horizontal axis) and $\log_{10}(p)$ (vertical axis).  The top-left
plot shows 500 points drawn from the prior described in the text.  The
top-right plot shows 500 points from the posterior distribution given
that someone with your memories exists, assuming $V=0.1$.  The bottom
plots show the posterior distributions assuming $V=1$ and $V=10$.
Tick marks are spaced one unit apart, representing change in $f$ or
$p$ by a factor of 10.  The diagonal lines indicate where
$\log_{10}(f)\,+\,\log_{10}(p)\, =\, -\log_{10}(V)$.

\end{figure}

The remaining plots in the figure show samples of points from the
posterior for $\log(p)$ and $\log(f)$ when $V$ is $0.1$, $1$, and
$10$.  The lines shown are where
$\log_{10}(p)\,+\,\log_{10}(f)\,=\,-\log_{10}(V)$, indicating for each
$f$ the value of $p$ that maximizes the probability that someone with
your memories exists.\footnote{Some details:\ \ 
The effect of the factor $p$ in
equation~(\ref{eq-you-given-fp}) is to shift the mean of
$\log_{10}(p)$ by $1.25^2\log(10)$, with the distribution remaining
Gaussian with the same standard deviation.  The remaining factor of
$\exp(-fpV)$ is never greater than one, so rejection sampling can be
used to obtain the posterior sample.}  Larger values of $V$ correspond
to a greater potential for another species to develop and then interfere
with your existence, a potential that is modulated by same factor,
$p$, that controls the likelihood of your development.  Accordingly,
larger values of $V$ shift the posterior distribution towards smaller
values of $p$.  The posterior distribution of $f$ is also shifted
towards smaller values (more so for large $V$), since smaller values
of $f$ reduce the probability that another species will interfere.

We can now determine the effect of FNC on your uncertainty concerning
one factor that enters into $p$.  Let us write $p=p_0p_1$, where $p_1$
is a single relevant factor, such as the probability that
multi-cellular life will evolve from single-celled life, and $p_0$ is
the product of all other factors.  Suppose that $p_0$ and $p_1$ are
independent, and that your prior distribution for $\log_{10}(p_1)$ is
Gaussian with mean $\log_{10}(0.1)$ and standard deviation 0.2, giving
a 95\% central interval for $p_1$ of 0.041 to 0.247, and a mean for
$p_1$ of $0.111$.  (Your prior for $\log_{10}(p_0)$ will therefore be
Gaussian with standard deviation $\sqrt{1.25^2-0.2^2}\,=\,1.234$.)
The conditional distribution for $\log_{10}(p_1)$ given $\log_{10}(p)$
is Gaussian with mean given by 
\beq
  \lefteqn{E[\log_{10}(p_1)\ |\ \log_{10}(p)]}\ \ \ \ \ \ \ \nonumber\\[4pt]
  & = &
    E[\log_{10}(p_1)]\ +\
    \big(\log_{10}(p)-E[\log_{10}(p)]\big)
      \,\Var[\log_{10}(p_1)]\,\big/\,\Var[\log_{10}(p)] \\[4pt]
  & = & \log_{10}(0.1)\,+\,\big(\log_{10}(p)-E[\log_{10}(p)]\big)\,0.2^2/1.25^2
    \\[4pt]
  & = & \log_{10}(0.1)\,+\,0.0256\,\big(\log_{10}(p)-E[\log_{10}(p)]\big)
  \label{eq-cond-m}
\eeq
and variance given by 
\beq
  \Var[\log_{10}(p_1)\ |\ \log_{10}(p)] & = &
  \Var[\log_{10}(p_1)]\,
    \big(1\,-\Var[\log_{10}(p_1)]\,\big/\,\Var[\log_{10}(p)]\big) \\[4pt]
  & = & 0.2^2\,(1\, -\, 0.2^2/1.25^2) \ \ =\ \ 0.0390
  \label{eq-cond-v}
\eeq
To find the posterior mean of $\log_{10}(p_1)$ given
that someone with your memories exists, we take the mean of~(\ref{eq-cond-m})
with respect to the posterior distribution of $\log_{10}(p)$.  The
posterior variance of $\log_{10}(p_1)$ is the sum of~(\ref{eq-cond-v})
and the variance of~(\ref{eq-cond-m})
with respect to the posterior distribution of $\log_{10}(p)$.

If $V=0$, so other intelligent species have no effect on earth, the
result of this computation is that the posterior mean and standard
deviation of $\log_{10}(p_1)$ are $\log_{10}(0.1236)$ and $0.2$, which
give a 95\% central interval of 0.050 to 0.305, and a posterior mean
for $p_1$ of $0.137$.  When $V=0$, the posterior distribution
of $\log_{10}(p_1)$ is Gaussian, and is the same as would be
obtained if $p_1$ were the only uncertain factor.  There is a
significant ``Presumptuous Philosopher'' effect from applying FNC,
although it is not as large in magnitude as some previous examples.

In contrast, the effect of FNC on the distribution of $p_1$ is much
less when $V$ is of significant size, even though, as can be seen in
the plots, the posterior distribution of $p$ itself is quite different
from the prior.  The posterior mean and standard deviation of
$\log_{10}(p_1)$ are $\log_{10}(0.1080)$ and $0.1985$ when $V=0.1$,
$\log_{10}(0.1042)$ and $0.1984$ when $V=1$, and $\log_{10}(0.1003)$
and $0.1984$ when $V=10$.  The posterior means of $p_1$ for these
values of $V$ are $0.120$, $0.116$, and $0.111$.  The last is nearly
identical to the prior mean of $p_1$, so there is almost no
``Presumptuous Philosopher'' effect on inference regarding this single
factor of $p$ when $V=10$.  With larger values of $V$, it is possible
for the posterior mean of $p_1$ to be less than the prior mean.

When $V=0$, the posterior distribution of $f$ is the same as the
prior, but with larger $V$, the posterior favours smaller values for
$f$, as can be seen in the plots.  We can look at a single factor
entering into $f$, just as we did for $p$.  If we write $f=f_0f_1$, we
can proceed much as above.  Suppose the prior for $\log_{10}(f_1)$ is
Gaussian with mean $\log_{10}(0.1)$ and standard deviation $0.2$,
giving a prior mean for $f_1$ of 0.111.  The posterior mean of $f_1$
is 0.111, 0.097, 0.093, and 0.089 for $V=0$, $V=0.1$, $V=1$, and
$V=10$.  A substantial change occurs when you condition on someone
with your memories existing, with the effect increasing as $V$
increases.

This is disturbing, since many of the factors contributing to $f$ ---
such as the probability of a technological civilization avoiding
self-destruction, and the probability that interstellar travel is
feasible --- are also relevant to human prospects, with larger values
being more favourable.  (However, some other factors going into $f$,
such as the probability that an intelligent species will decide to
destroy the potential habitat of another intelligent species, are 
ones that many of us would not wish to be large.)  There is thus a
``doomsday'' aspect to this analysis, since use of FNC has revealed
that we should increase the probability we assign to some negative
scenarios, above the probability we would assign based on ordinary
considerations.  The source of this pessimism is quite different from
that of the Doomsday Argument of Section~\ref{sec-da}, however.  It is
based on the empirical observation that we are not aware of any other
intelligent species.  One possible explanation of this observation is
that most intelligent species are destroyed in some fashion, or at
least fail to develop in a way that would make their presence known to
us.  This is a reason to increase our assessment of the probability of
this happening to us.

The magnitude of pessimism that this argument warrants depends on our
beliefs regarding a wide range of topics in physics, astronomy,
biology, and sociology.  In contrast, the Doomsday Argument depends
only on the size of the future human population in different
scenarios, and can produce very large probabilities of ``doom'' if the
alternative is a future involving interstellar colonization, or even
just intensive settlement of the solar system.  Arguments based on FNC
are unlikely to produce such extreme pessimism.

\subsection{\hspace*{-7pt}Inflation and infinite
                          universes}\label{sec-inf}\vspace*{-5pt}

Cosmological theories in which an early period of ``inflation''
greatly expands the universe imply that the universe we are in is
infinite in size.\footnote{At least according to Knobe, Olum, and
Vilenkin (2006), though Olum (2004) says models of finite inflationary
universes can be contrived.}  Furthermore, in most such cosmologies,
our universe is only one of many within a larger ``multiverse''.
Finally, these theories do not produce any tight linkage between
distant parts of the universe, which might constrain them either to be
similar in detail, or different.  It follows that in an inflationary
universe we should expect all possible observers to exist, each an
infinite number of times.  In particular, there should be an infinite
number of distant observers with exactly your current memories.

This is a problem for FNC.  If you accept inflationary cosmology as
correct, someone with your memories exists with probability one,
regardless of what else might be true.  Conditioning on the existence
of someone with your memories will then have no effect on the
probabilities of any other theories.  In particular, FNC no longer
provides a solution to the Freak Observers problem
(Section~\ref{sec-fo}).

However, an infinite universe leads to many other problems as well.
For example, Knobe, Olum, and Vilenkin (2006) discuss the ethical
implications of an infinite universe.  These and other issues with
infinities do not seem at all clear to me.  As an interim solution, I
advocate simply \textit{ignoring} the problem of infinity, which is
certainly what everyone does in everyday life.

My primary justification for ignoring the problems FNC has with
infinity is that the finite size necessary to cause problems is
extraordinary large.  As discussed in Section~\ref{sec-fo}, for
duplicate observers to exist, the universe must be a huge number of
orders of magnitude greater than the portion of the universe we
presently observe.  Can the difference between an unimaginably vast
universe and one that is truly infinite actually be crucial to our
inferences regarding local matters, which concern only the small
region within a few tens of billions of light years of us?

Conceivably, the answer is yes.  But it seems more likely to me that
either the infinity will disappear once the theory is better
understood, or it will turn out that its implications, at least for
the questions dealt with here, are not great.  The situation is
analogous to thought experiments with extreme assumptions, where (as
discussed in Section~\ref{sec-fant}) there is a danger that our
reasoning will implicitly use premises that are not true given these
extreme assumptions.  The difference, of course, is that the extreme
assumption in this case concerns reality, and may ultimately prove
unavoidable.  But it seems best to try to avoid it at least initially.

Some technical matters also support the strategy of ignoring infinity.
First, even if the universe we are in is infinite, our knowledge of it
is certainly not infinite, since distant parts of the universe are
outside our past light-cone, due to the universe's finite age.  This
is fortunate, since if we were subjected to non-negligible influences
from every part of an infinite universe, our experience would be a
incomprehensible jumble.  (This is just a more general form of Olbers'
Paradox --- that if the universe is infinite, the night sky should be
white.)  Should the infinite regions with which we have had no contact
really count as part of ``our'' universe?  One might argue that they
should, one the grounds the we might be in contact with them in the
future.  Whether this is so depends on details of the universe's
expansion, but let's suppose that any two regions of the universe,
even very distant ones, will eventually come into contact.  Who will
receive information from such distant regions?  You will likely be
dead, but suppose instead that you have achieved immortality.  If you
are actually attending to news from distant regions, you must be
expanding your memory.  But any increase in your memory results in a
huge increase in the size of universe needed for a duplicate observer
with your exact memories to exist.  So it is difficult to imagine any
scenario in which the existence of duplicate observers has
observational consequences.

It is therefore not surprising that the puzzle presented by Olum
(2004) as arising from inflationary cosmology is not really dependent
on the universe being infinite.\footnote{A universe of the size we
observe probably suffices.  A bigger universe could be necessary if
for some reason life is extraordinary rare, but if so, the larger size
will not cause problems for FNC, since the probability of duplicate
observers will also be lower if life is rare.  It is possible that
Olum sees an infinite universe as necessary to justify use of SSA,
thinking it would then be the only way to avoid the Freak Observers
problem.}  Olum considers the probability that an intelligent species
will colonizing its galaxy (or even many galaxies), and thereby
achieve an enormous population (eg, $10^{19}$ individuals), and
concludes that the probability of a species doing this, while perhaps
substantially less than one, is not minuscule.  Accordingly, one would
expect most intelligent observers to belong to such a galactic
civilization.  Yet we don't.  Olum sees this as a conflict between
observation and ``anthropic reasoning'', by which he means
SSA$-$SIA.\footnote{Previously (Olum 2002), he had advocated SSA+SIA,
but he apparently had doubts about SIA at the time of this paper, and
more recently (Knobe, Olum, and Vilenkin 2006).}  If you consider
yourself to be a randomly selected observer, as advocated by SSA, you
should very likely be either a member of some other species that has
colonized their galaxy, or be a human from later in our history, when
we have done so.

Although Olum doesn't present it as such, this is essentially the
Doomsday Argument (applied in the context of many species, as in
Section~\ref{sec-da-et}), except that Olum is sufficiently confident
that doom is not nearly universal that he regards the result of the
argument as a paradox rather than a prediction.  Trying to resolve the
paradox, he considers the possibility that ``anthropic reasoning''
(SSA$-$SIA) is invalid, along with other possibilities (eg, galactic
colonization is actually exceedingly difficult).  He comes to no
definite conclusion, but considers that several of these possibilities
might together be sufficient to explain the conflict.

My conclusion is that ``anthropic reasoning'' --- meaning SSA$-$SIA
--- is indeed invalid.  In contrast, application of FNC produces no
paradox.  Suppose we know that intelligent species are very rare, and
hence seldom interact, so that we can ignore the complexities
discussed in Section~\ref{sec-int}.  The probability of an observer
such as you existing, as a member of a species that has not colonized
the galaxy, is determined by factors influencing the evolution of
species up to our stage of development.  It is irrelevant what happens
to these species later; hence it is irrelevant whether, for example,
galactic colonization is easy or hard.\footnote{Two explanations
considered by Olum are not irrelevant:\ \ We are part of a galactic
civilization without being aware of it, or (perhaps a special case of
this) we don't actually live on earth, but rather exist in a computer
simulation.  Both possibilities could increase the probability of
someone with your memories existing, and hence might be favoured by
FNC.}

Anthropic reasoning has also been applied to theories in which
multiple inflating universes can have different values of fundamental
physical parameters, in particular the ``cosmological constant'',
which influences how rapidly the universe expands.  The observed value
of the cosmological constant is close to, but not exactly, zero.  The
most well-accepted theories of the cosmological constant provide no
apparent reason for it to be as small as it is --- it might equally
well have any value over a range that is $10^{120}$ larger than its
actual value.  However, Weinberg (1987) calculated that only a much
narrower range of values around zero will lead to the formation of
galaxies, which he considered a prerequisite for life to exist.  Since
subsequent measurements found a non-zero value in this range, this
calculation has been taken to be a successful prediction using the
Anthropic Principle.  I will critique such reasoning below in
connection with string theory, where this and related anthropic
arguments have recently become prominent.

\subsection{\hspace*{-7pt}The landscape of string
                          theory}\label{sec-land}\vspace*{-5pt}

String theory is an attempt to unify Einstein's theory of gravity with
the ``Standard Model'', which describes electromagnetism and the weak
and strong nuclear forces.  String theorists originally hoped that the
requirement of mathematical consistency would produce a unique theory,
which would predict the previously unexplained parameters of the
Standard Model, such as the masses of elementary particles.  Though
this possibility has not been definitely ruled out, many string
theorists now think it more likely that hundreds of parameters of the
theory can be varied while retaining consistency.  This results in a
huge ``landscape'' of possible physical laws, with perhaps $10^{500}$
or more possibilities, each of which produces different values for the
parameters of the Standard Model, and for the cosmological constant.
The universe sits in a ``valley'' in this landscape, to which it
``descended'' during its inflationary period.  If a huge number of
inflating universes formed, or if, following the Many Worlds
interpretation of quantum mechanics, a single universe has a huge
number of superposed versions, almost all valleys of the landscape
will be populated by one or more universes.  The landscape will then
describe not just a set of mathematical possibilities, but an actual
multitude of real universes.

This view of string theory and cosmology has been advocated by
Susskind (2003, 2006), who then uses it as a basis for anthropic
reasoning.  In his view, the particular values of the physical
parameters we observe, and indeed, even the set of particles we see,
cannot be explained by the requirements of any fundamental theory, but
they can be explained (at least partially) by the requirement that the
universe contain intelligent observers such as ourselves, since
otherwise there would be no one to look for an explanation.  In other
words, we measure the values of the physical constants that we do
because only these values (or similar values) allow for the existence
of someone to measure them.

Unfortunately, Susskind is not too clear on the exact purpose of this
reasoning, or its justification.  Before attempting to critique his
views, I will try to clarify the issues by discussing what one might
conclude by applying FNC.

As discussed in Section~\ref{sec-inf}, infinite universes pose a
problem for FNC.  Accordingly, I will suppose that the landscape is
populated by only around $10^{600}$ universes (more than enough to
guarantee at least one in each of $10^{500}$ valleys) and that each
universe has at most $10^{350}$ galaxies (much more than the $10^{11}$
we can observe in our universe).  If each galaxy has $10^{10}$
inhabited planets, each of which has a generous $10^{20}$ inhabitants,
who are replaced by equal numbers for $10^{20}$ generations, the total
number of intelligent observers in all universes who ever exist will
be at most $10^{1000}$.  As discussed in Section~\ref{sec-fo}, this is
far too few for there to be any non-negligible chance of another
observer with your exact memories existing.  

In this scenario, we can apply FNC without difficulty.  Conceivably,
the answers we obtain might not be correct if in reality there are an
infinite number of universes of truly infinite extent.  However, in
none of the anthropic arguments I am aware of does such a distinction
between unimaginably vast and truly infinite universes play any
apparent role.  If infinity is actually crucial, the proponents of
anthropic arguments need to make the reason for this more explicit.

Consider two cosmological theories, in both of which there are
$10^{600}$ universes formed by inflation.  In theory $L$, these
universes populate a huge number of valleys in a landscape of possible
physical laws, as described above.  In theory $S$, there is either no
landscape, perhaps because the requirement of mathematical consistency
uniquely determines physical laws, or the landscape has only a single
valley, which has much the same effect once inflation is over.  If
string theory is accepted as the correct foundation of physics, and
its basic principles are not in dispute, whether $L$ or $S$ is the
correct theory may be a mathematical question, whose answer we are
ignorant of only because of the difficulty of performing the necessary
calculations.  Alternatively, $L$ and $S$ may have different
foundations, even though they both lead to similar collections of
inflating universes.  In either case, suppose that, on mathematical
and physical grounds, you judge the two theories to be equally
plausible.  What should you judge the probabilities of these theories
to be after applying FNC, conditioning on all your memories, both of
everyday life, and of the results of whatever scientific measurements
have been performed?

We can distinguish two situations.  First, suppose that the unique
parameters underlying theory $S$ are known, and that at least some of
the implications of these fundamental parameters for the parameters of
the Standard Model and for the value of the cosmological constant have
been worked out.  I'll refer to this version of theory $S$ as theory
$S_1$.  If the implications of theory $S_1$ contradict experimental
measurements, we can clearly rule out $S_1$, and conclude that
theory $L$ is true (assuming that these are the only alternatives).
Note that ``measurement'' includes general observations, such as the
existence of galaxies, which may rule out certain values for
parameters of the Standard Model or for the cosmological constant.
Alternatively, the fundamental parameters of theory $S_1$ may produced
parameters for the Standard Model and cosmological constant that lie
within the region, $Y$, that so far as you know is compatible with
observation.  The probability that someone with your memories will
exist according to theory $S_1$ will then be $10^{600}$ times the
probability that a universe with parameters in $Y$ will produce an
observer with your memories.\footnote{This probability is (at least
roughly) the same for all universes with parameters in $Y$, since $Y$
is defined to be the region of parameters that can't be ruled out
based on your memories.}  On the other hand, the probability that
someone exists with your memories under theory $L$ will be $10^{600}$
times the fraction of valleys in the landscape that produce parameters
in region $Y$ times the probability that a universe with parameters in
$Y$ will produce an observer with your memories.  The odds in favour
of theory $L$ will therefore be equal to the fraction of valleys in
the landscape that produce parameters in $Y$.\footnote{It's possible
that theory $L$ defines some non-uniform measure over valleys, in
which case the odds would be the total measure for valleys in region
$Y$ divided by the total measure for all valleys.  This elaboration
does not affect the basic argument.}  The landscape of string theory
is typically seen as containing valleys with a great diversity of
physical laws, so the odds in favour of $L$ in this scenario would be
tiny --- ie, theory $S_1$ would be very strongly favoured.

In the second situation, the fundamental parameters for theory $S$ are
unknown.  Perhaps, for example, it has been proved that the structure
of theory $S$ (but not theory $L$) must lead to a unique set of
parameters, but their actual values are not known, though mathematical
intuition allows one to give them some broad prior distribution.  Or
theory $S$ might just baldly state that the universes that exist have
only a single set of parameter values, but these values are arbitrary,
with some broad prior distribution.  I'll use $S_*$ to refer to a
theory $S$ of this type.  In this situation, the probability that
someone with your memories exists under $S_*$ is equal to $10^{600}$
times the prior probability of region $Y$ times the probability that
someone with your memories will exist in a universe whose parameters
are in $Y$.  The odds in favour of theory $L$ will be equal to the
fraction of valleys that produce parameters in $Y$ according to $L$
divided by the prior probability of $Y$ according to $S$.  If the
distribution of parameters of valleys in $L$ is similar in breadth to
the prior for parameters in $S_*$, these odds will be roughly one ---
ie, neither theory will be favoured, since neither gives any very
specific predictions.

In these arguments, a crucial role is played by the region $Y$, which
encompasses values of the parameters of the Standard Model and of the
cosmological constant that are not ruled out by your memories
(including your memories of scientific measurements).  In contrast, it
is irrelevant what region of parameters is compatible with life, or
with intelligent life, or with intelligent life that has developed a
scientific culture.  These regions would likely be much bigger than
$Y$, since there is no apparent reason why, for instance, life
couldn't develop in a universe with only half as many galaxies as we
see.

These applications of FNC accord with usual scientific reasoning.  If
theory $S$ makes specific predictions, and these are compatible with
what is observed, it is favoured over theory $L$, since $L$ makes no
specific predictions.  If theory $S$ also makes no specific
predictions, either because it has not been sufficiently worked out,
or because it has arbitrary parameters, then neither $S$ nor $L$ are
favoured.

How is this different from ``anthropic'' reasoning?  The crucial point
seems to be that theory $S_1$, whose parameters are known, and match
observations, implies that $10^{600}$ universes much like ours exist.
In contrast, theory $L$ implies that far fewer universes like ours, or
even compatible with life, will exist.  (Though it is assumed that
theory $L$ implies the existence of at least one universe with
intelligent life.)  If one applied SSA+SIA, the prior probability of
theory $S_1$ would be greatly boosted compared to that of theory $L$,
and the result would be the same as applying FNC.  But if one instead
applies SSA$-$SIA, there is no boost for theory $S_1$.  The crux of
the ``anthropic'' argument seems to be that one should not penalize
theory $L$ for predicting that only a few habitable universes exist,
as long as it predicts at least one, since we will naturally find
ourselves in a habitable universe, even if they are rare.  As a
result, the degree to which theory $S_1$ is favoured over theory $L$
is much reduced.

A numerical example may clarify the situation.  Suppose that the
landscape of theory $L$ has $10^{500}$ valleys, whereas theory $S_1$
has only one valley, whose properties are known.  The single valley of
theory $S_1$ is compatible with intelligent life, and furthermore,
with your specific memories.  Of the $10^{500}$ valleys of theory $L$,
$10^{10}$ are compatible with intelligent life, and $10^{6}$ of these
are compatible with your specific memories.  For simplicity, let's
assume that all universes with intelligent life have the same
population.  Application of FNC then gives odds of
$10^{6}/10^{500}\,=\,10^{-494}$ for theory $L$, but the anthropic
reasoning described above, based on SSA$-$SIA, gives odds of
$10^{6}/10^{10}\,=\,10^{-4}$ for $L$.  So whereas $L$ is essentially
disproved if FNC is used, it retains a non-negligible probability
under SSA$-$SIA.  This result may seem reasonable if you take an
anthropic view, but note the disturbing sensitivity of the odds for
$L$ to the definition of ``intelligent life'', and the need to
determine whether such life exists in universes with $10^{500}$
different physical laws before a conclusion can be reached.

This situation resembles that discussed in Sections~\ref{sec-density}
and~\ref{sec-density-fnc}, where Marochnik's theory and other theories
in which the density of observers vary were considered.  I argued
there that the results of FNC and SSA+SIA are correct using the device
of companion observers.  In the context of inflationary cosmology, the
``companions'' would need to exist in every universe, regardless of
whether it is hospitable to us, even though the physical laws differ
radically from universe to universe.  Suppose there are a great many
such observers in every universe (albeit unobserved by us, so far),
and that they know that observers like them exist in every universe.
They will take the existence of humans in this universe as evidence
that many universes have physical laws that allow beings like humans
to exist.  Why should we think differently?  One might well wonder
whether this scenario is stretching the concept of companion observers
too far, but I see no specific reason for thinking that these
conclusions are inappropriate.

Susskind does not discuss anthropic reasoning in terms of probabilistic
principles such as SSA$-$SIA, nor in reference to any explicit comparison 
of theories.  Instead, his focus seems to be on
finding an \textit{explanation} for our universe's physical laws.
In his book, \textit{The Cosmic Landscape}, he describes how he
came to accept use of the Anthropic Principle, beginning with an
account of the many ``coincidences'' that seem necessary for life to 
exist:
\bqt\indent
   There are multiple ways that things could go wrong with the
   nuclear cooking [~of heavy elements~]. \ldots{} But again, it
   would do no good for the nuclear physics to be ``just 
   right'' if the universe had no stars. \ldots{} How then
   did the universe get to have such a large preponderance of
   matter over antimatter? \ldots{} Another essential requirement
   for life is that gravity be extremely weak. \ldots

   Just how seriously should we take this collection of 
   lucky coincidences?  Do they really make a strong case for
   some kind of Anthropic Principle?  My own feeling is that they 
   are very impressive, but no so impressive that they would
   have pushed me past the tipping point to embrace an
   anthropic explanation. \ldots{} accidents, after all, do happen.

   However, the smallness of the cosmological constant is another
   matter.  To make the first 119 decimal places of the vacuum
   energy zero is most certainly no accident.  But it was not
   just that the cosmological constant was very small.  Had
   it been even smaller than it is, had it continued to be zero
   to the current level of accuracy, one could have gone on
   believing that some unknown mathematical principle would
   make it exactly zero. \ldots

   But even the cosmological constant would not have been enough
   to tip the balance for me.  For me the tipping point came
   with the discovery of the huge Landscape that String Theory
   appears to be forcing on us.  (Susskind 2006, pp.~182-185)
\eqt
I will discuss the cosmological constant in more detail below, but
for now let us count it as just one more ``lucky coincidence''.
The last point above seems crucial.  He expands on it later:
\bqt
   \ldots{} in my own mind, the ``straw that broke the camel's
   back'' was the realization that String Theory was moving
   in what seemed to be a perverse direction.  Instead of
   zeroing in on a single, unique system of physical laws,
   it was yielding an ever-expanding collection of Rube 
   Goldberg concoctions.  I felt that the goal of a single
   unique string world was an ever-receding mirage and
   that the theorists looking for such a unique world were on
   a doomed mission.  (Susskind 2006, p.~199)
\eqt

In terms of theories of type $S$ and $L$ discussed above, it appears
that Susskind initially saw string theory as a theory $S_*$, for which
it was believed (though not proven) that the fundamental parameters of
the theory were unique, even though they were unknown.  If he had
thought to compare it to some theory $L$ (obviously based on some
different structure), and had applied FNC, he would have concluded
that the two theories were about equally likely, since at that point
neither could make specific predictions.  Of course, he would have
hoped to find the unknown unique parameters of $S$, and he would have
hoped that the predictions of theory $S$ with these parameters matched
observations.  If both hopes had been fulfilled, application of FNC
would have produced the conclusion that this theory (now of type $S_1$)
was vastly more probable than theory $L$.  Perhaps these applications
of FNC approximate the logic Susskind employed at that time.

After abandoning the quest for a unique set of physical laws, accepting 
instead a landscape of possible laws, populated by multiple universes,
Susskind appears to have been concerned with only one competing theory
--- that the particular laws we see were chosen by an Intelligent
Designer, with the purpose of creating a universe containing life.  It
is this alternative that his anthropic arguments appear aimed at
refuting, or at least rendering unnecessary.  In his book, which is
subtitled ``The Illusion of Intelligent Design'', he writes: 
\bqt\indent
   To Victor's [~a friend's~] question, ``Was it not God's
   infinite kindness and love that permitted our existence?''
   I would have to answer with Laplace's reply to Napol\'{e}on:
   ``I have no need of this hypothesis.''\ \ \textit{The Cosmic
   Landscape} is my \mbox{answer}\ldots{} (Susskind 2006, p.~15)
\eqt

Obtaining this answer doesn't require anthropic reasoning, however.
Intelligent Design can be seen as a theory $S_*$ in which all
universes operate by a single set of physical laws that are fixed to
arbitrarily values by the Designer.  Supposing we have some broad
prior distribution for the parameters of these physical laws, we find
that the theory makes no specific predictions.  Application of FNC
leads to the conclusion that this theory and theory $L$ are about
equally likely.  There is ``no need'' for the hypothesis of an
Intelligent Designer.

An advocate of Intelligent Design might, of course, maintain that a
broad prior is not appropriate --- that the prior should be confined
to physical laws that will produce a universe containing intelligent
life.  I'll call this theory $S_D$.  If you consider $S_D$ and $L$
equally likely \textit{a~priori}, FNC will lead you to conclude that
theory $S_D$ is much more probable than theory $L$ --- extending the
numerical example above, theory $S_D$ will predict that the universe
has one of the $10^{10}$ sets of laws that are compatible with
intelligent life, of which $10^{6}$ are compatible with your
observations, so the odds in favour of $L$ will be
$(10^{6}/10^{500})\,/\,(10^{6}/10^{10})\ =\ 10^{-490}$.  But why
should one think the Designer wished intelligent life to exist, as one
must to regard $S_D$ as plausible?  Some may think this, but an
argument that has as a premise God's infinite kindness and love for
humanity is not a scientific argument, and requires no scientific
refutation.

Nevertheless, if one wishes a counter-argument, anthropic reasoning
may appear to provide one.  Applying SSA$-$SIA will make the theories
of the landscape ($L$) and of an Intelligent Designer who likes
intelligent life ($S_D$) equally likely.  Of the $10^{500}$ valleys in
theory $L$, only the $10^{10}$ with intelligent life ``count'' when
applying SSA$-$SIA, so the probability of a universe compatible with
what you observe is $10^{6}/10^{10}$, the same as for theory $S_D$.

I have argued in this paper that SSA$-$SIA is not a valid principle of
reasoning.  If so, one would expect Susskind's approach to produce
strange results in other contexts.  Consider a comparison of
Susskind's theory $L$, in which there is a landscape of $10^{500}$
possible physical laws, with a theory $S_*$ that other string
theorists may still be working on, in which it is thought that only
one of these $10^{500}$ apparent possibilities is mathematically
consistent, though it is not known which of the $10^{500}$ it is.  As
discussed above, applying FNC leads to the conclusion that $S_*$ and
$L$ are equally likely.  What is the result of applying SSA$-$SIA?

SSA$-$SIA will strongly favour theory $L$.  In the numerical example
above, theory $L$ would predict a universe compatible with what you
see with probability $10^6/10^{10}\, =\, 10^{-4}$, since of the
$10^{10}$ valleys in the landscape that allow intelligent life, $10^6$
are compatible with your observations.  The corresponding probability
under theory $S_*$ is only $10^6/10^{500}\, =\, 10^{-494}$, so it is
very strongly disfavoured, with the odds in favour of $L$ being
$10^{-4}/10^{-494}\,=\,10^{490}$.  Another way of looking at this
problem is to split theory $S_*$ into theories
$S^1,\,S^2,\,\ldots,\,S^{10^{500}}$, one for each possible set of
physical laws, and split the prior probability of $1/2$ for $S_*$ into
prior probabilities of $0.5\times10^{-500}$ for each of these
theories.  All but $10^6$ of these theories are incompatible with your
observations.  The total posterior probability of all the sub-theories
of $S_*$ that are compatible with what you see works out to
$10^6\times0.5\times10^{-500}\ /\ (10^6\times0.5\times10^{-500}\ +\
0.5\times10^{-4})\ \approx\ 10^{-490}$.

This seems unreasonable.  Perhaps there are good reasons to think that
the old research programme of looking for unique physical laws within
string theory has poor prospects, but until it is actually proved
hopeless, its chances of success are surely not as low as $10^{-490}$.
Susskind does not explicitly draw such a pessimistic conclusion, but
it seems to follow from the logic of anthropic reasoning that he uses.

Susskind does draw an even more surprising conclusion from the
anthropic viewpoint.  Discussing the idea that the laws of physics
might be an emergent phenomenon, of the sort that is well-known in
condensed-matter physics, he writes: 
\bqt
  The properties of emergent systems are not very flexible.
  There may be an enormous variety of starting points for the
  microscopic behavior of atoms, but\ldots{} they tend to lead
  to a very small number of large-scale endpoints. \ldots{}
  This insensitivity to the microscopic starting point is the
  thing that condensed-matter physicists like best about 
  emergent systems.  But the probability that out of the small
  number of possible fixed points (endpoints) there should be
  one with the incredibly fine-tuned properties of our anthropic
  world is negligible. \ldots{} A universe based on conventional
  condensed-matter emergence seems to me to be a dead-end idea.
  (Susskind 2006, pp.~359-360) 
\eqt 
This comment is remarkable.  An inflexible theory leading to only a
small number of possible sets of physical laws (preferably just one)
is what Susskind had originally hoped string theory would be!  Yet now
he sees such a theory as being almost certainly false, not (just)
because of detailed problems with it, but because of the very
inflexibility, leading to near uniqueness, that he previously saw as
one of the most attractive features of string theory.  Moreover, even
application of SSA$-$SIA does not lead to this theory being greatly
disfavoured, if the details have not been worked out that would show
what the small number of possibilities actually predict. 
Rather than $10^{500}$ sub-theories as in the example above, there are only, 
say, 10 sub-theories, each of which has a substantial portion of the
prior probability for the theory as a whole.
The low probability Susskind assigns to this theory can only come from
his assigning a low \textit{prior} probability to the whole theory,
based on a prior belief that physical laws do not
have any simple explanation, but are instead a ``Rube Goldberg concoction''.

Such a belief is, of course, contradicted by numerous scientific
success stories, such as the use of quantum mechanics to explain the
complex features of atomic spectra.  However, some other complex
phenomena do seem to have no explanation other than accident --- the
outlines of the continents, for example, have no fundamental
geological explanation.  Whether a phenomenon has a simple explanation
or not cannot be determined \textit{a~priori}.  Perhaps a multiplicity
of universes with differing physical laws exist; perhaps the set of
possible physical laws is much more constrained.  One can tell which
only by creating and testing theories of both sorts.

Anthropic reasoning has also been criticized by Smolin (2006), who has
in addition proposed a third possibility --- universes with a great
diversity of physical laws can indeed exist, but rather than the
physical laws of each universe being chosen at random from some simple
distribution, they are chosen according to some dynamical process,
which leads to a distribution of universes in which the physical laws
we actually observe are much more likely.  He proposes a particular
theory of ``cosmic natural selection'', based on the idea that new
universes are formed inside black holes, with slightly perturbed
physical laws.  Selection will then tend to favour physical laws that
make a universe produce many black holes.  Smolin argues that such a
universe will resemble ours.  

A successful theory of this sort would be greatly favoured by FNC, in
comparison with a theory that distributes universes uniformly over
valleys of the landscape, since it would (if successful) greatly
increase the probability of a universe similar to ours (in the region
$Y$ defined earlier), and hence also the probability that someone with
your memories will exist.  In contrast, such a theory might not be
favoured at all by SSA$-$SIA.  Universes without intelligent life
``don't count'' with SSA$-$SIA, so if Smolin's theory (for example)
leads to many more universes that contain intelligent life, but fails
to further concentrate the distribution towards universes more
precisely like ours, it will not be considered more probable by
SSA$-$SIA than a theory in which the physical laws for each universe
are drawn from a much broader distribution.

As promised above, I will now consider in more detail the issue of the
cosmological constant, which is usually denoted by $\Lambda$.  As seen
from the quote above, Susskind considers the observed small, but
non-zero, value for $\Lambda$ to be the strongest of the
``coincidences'' that led him to consider anthropic explanations.  Two
separate aspects of the situation seem responsible for this --- the
large magnitude of the coincidence, and the special role of the value
zero.

The range of values for $\Lambda$ that are compatible with life (taken
to be the range for which galaxies form) is much narrower than the
range of values that seem plausible on general theoretical grounds, by
a factor of roughly $10^{120}$.  This ratio of prior range to
``anthropic'' range (for which life can exist) is substantially
greater than for the other parameters of the Standard Model that seem
to be ``fine-tuned'' for life.  Someone who accepts the basic
anthropic argument (based, so far as I can tell, on SSA$-$SIA) will
naturally be impressed by this.  As I argue above, however,
application of FNC does not lead one to favour a theory based on the
landscape for this reason, so the magnitude of the coincidence is
irrelevant.

Does the fact that the observed value of $\Lambda$ is close to zero,
but not exactly zero, modify this conclusion?  Consider some other
parameter, for which the range of conceivable values is $(0,1)$ and
the range of values compatible with life is $(0.3181,0.3192)$.  The
best measurement of this parameter gives the 95\% confidence interval
$(0.3185,0.3187)$.  Suppose you consider an anthropic explanation for
the value of this parameter to be attractive.  Someone now advances a
plausible theory that the true value is exactly $1/\pi=0.3183\ldots$,
which is somewhat at variance with the measurement, but not hopelessly
so.  After learning of this theory, should your confidence in an
anthropic explanation be greater or less than before?  Surely you
should be less confident, since it's possible that this new theory
provides the true explanation.  Certainly you shouldn't be
\textit{more} confident in an anthropic explanation now that before.

Analogously, the fact that the anthropic range for the cosmological
constant includes the special value zero, which one might imagine
could result from some theoretical constraint enforcing cancellation
of terms, does \textit{not} make an anthropic explanation more likely,
but rather the reverse.  This is partly because of the possibility
that $\Lambda$ is indeed exactly zero, even though current
observations indicate otherwise.  More likely, however, is that some
theory might explain why $\Lambda$ is close to, but not exactly, zero.

Even if no good non-anthropic explanation for $\Lambda$ being near
zero can be found, the anthropic explanation may have its own problem
--- why is the special value zero contained in the rather narrow
anthropic range (about $10^2$ wide compared with a prior range of
$10^{120}$)?  The anthropic range of $\Lambda$ is a function of the
\textit{other} parameters of the physical laws.  Why should these
other parameters conspire to make this range contain zero?  Perhaps
there is some plausible cosmological explanation, valid even when the
set of particles is much different from what we observe, but I have
not seen the issue discussed.

To summarize, at least the following seem possible explanations for the 
value of $\Lambda$:\vspace{-6pt}
\begin{itemize}
\item[1)] $\Lambda$ must be exactly zero, for theoretical reasons.
\item[2)] $\Lambda$ has a specific value that is close to but not equal to
        zero, for theoretical reasons.
\item[3)] $\Lambda$ has a value that is not completely determined theoretically,
        but which theory says is likely to be close to zero.
\item[4)] $\Lambda$ takes on various values in different valleys of
        the landscape.  A non-negligible fraction of these
        values are close to zero, with the others being widely distributed.
\item[5)] $\Lambda$ takes on various values in different valleys of
        the landscape, with no tendency for these values to be close to 
        zero.\vspace{-6pt}
\end{itemize}
Explanation~(1) is viable only if current observations are in error.
A theory of the sort required for explanation~(2) would seem on
general grounds to be conceivable --- reasons for something to be zero
often can be modified to produced reasons for something to be near
zero.  Explanations~(3) and~(4) are not entirely distinct.  Smolin's
cosmic natural selection theory (Smolin 2006) and a recent cyclic
model of the universe due to Steinhardt and Turok (2006) provide
explanations of this type.  Explanation~(5) provides a reason to think
that a small value of $\Lambda$ is possible, but explains why we see
such a rare value only if you accept anthropic explanations.

Arguments in favour of an anthropic explanation for the cosmological
constant seem to generally dismiss explanations~(1) to~(4), though
Weinberg (2000) remarks that an \textit{a priori} distribution for
$\Lambda$ with a peak near zero would obviate anthropic explanations.
If only explanation~(5) is considered, however, anthropic reasoning
does no actual work, but just makes one feel more comfortable.  The
real question is whether the Anthropic Principle provides good reason
to increase the probability of explanation~(5) compared to the others.
The effect of applying SSA$-$SIA, as discussed earlier, is to let a
theory predict many lifeless universes (eg, with $\Lambda \gg 0$)
without penalty, as long as it predicts at least one universe with
intelligent life.  In contrast, when FNC is applied, explanation~(5)
is heavily penalized compared to an otherwise plausible theory that
provides an explanation of type~(1) to~(4).  (This assumes that all
theories produces a similar collection of universes --- if not, we get
into the difficult problem (discussed in Section~\ref{sec-varsize}) of
comparing theories that differ in the size or multiplicity of
universes.)

My conclusion is that when FNC can be clearly applied, it does not
support the type of anthropic reasoning that has been used to
``explain'' the apparent fine-tuning of physical constants to values
necessary for life, via a multitude of universes populating a
landscape of physical laws.  Such anthropic reasoning appears to be
based on SSA$-$SIA, and shares with it a disturbing sensitivity to the
reference class chosen.  Moreover, SSA$-$SIA, in both this application
and its applications to the problems discussed previously, can produce
counterintuitive conclusions.  When the universe is truly infinite,
and especially when different theories predict universes of different
sizes, some more general version of FNC is needed.  However, I see no
reason at present to think that my conclusions regarding anthropic
reasoning would be invalid in these situations.

None of this says that a cosmology with multiple universes populating
a landscape of physical laws cannot be correct.  FNC does give a
preference to theories of this sort in which the distribution of
universes is concentrated on valleys in the landscape that produce the
physical laws we see, but perhaps no such theory is viable.  Many
ordinary phenomena, such as the outlines of the continents and the
radii of the orbits of the planets, are believed to have no
explanation other than accident.  On cannot rule out \textit{a priori}
the possibility that the cosmological constant and the parameters of
the Standard Model have only this explanation, since this might be the
truth.  Such an explanation is, however, a ``last resort'', in that
any theory that more specifically predicts the observered values, and
is otherwise acceptable, should be greatly preferred.

\section*{Acknowledgements}\vspace{-10pt}

I thank David MacKay for helpful discussions.  This research was
supported by the Natural Sciences and Engineering Research Council of
Canada.  I hold a Canada Research Chair in Statistics and Machine
Learning.

\section*{References}\vspace*{-10pt}

\leftmargini 0.2in
\labelsep 0in

\begin{description}
\itemsep 2pt

\item[] Bostrom, N.\ (2002) \textit{Anthropic Bias:\ \ Observation
  Selection Effects in Science and Philosophy}, New York: Routledge.

\item[] Bostrom, N.\ (2006) ``Sleeping Beauty and self-location: A
  hybrid model'', preprint available from \texttt{http://www.nickbostrom.com},
  to appear in \textit{Synthese}.

\item[] Bostrom, N.\ and \'{C}irkovi\'{c}, M.~M.\ (2003) 
  ``The doomsday argument and the self-indication assumption: reply to
  Olum'', \textit{The Philosophical Quarterly},
  vol.~53, pp.~83-91.

\item[] Brin, D.\ (1983) ``The `Great Silence': The controversy 
  concerning extraterrestrial life'', \textit{Quarterly Journal of
  the Royal Astronomical Society}, vol.~24, pp.~283-309.

\item[] Carter, B.\ (1974) ``Large number coincidences and the 
  Anthropic Principle in cosmology'', in M.~S.~Longair (editor),
  \textit{Confrontation of Cosmological Theories with Data},
  Dordrecht: Reidel, pp.~291-298.

\item[] Carter, B.\ (2004) ``Anthropic principle in cosmology'', 
  available from the e-print archive at 
  \texttt{http://arxiv.org/abs/gr-qc/0606117}.

\item[] Elga, A.\ (2000) ``Self-locating belief and the Sleeping
  Beauty problem'', \textit{Analysis}, vol.~60, pp.~143-147.

\item[] Gott, J.~R.\ (1993) ``Implications of the Copernican Principle for 
  our future prospects'', \textit{Nature}, vol.~363, pp.~315-319.

\item[] Hanson, R.\ (1998) ``Burning the cosmic commons: Evolutionary
  strategies of interstellar colonization'', available from
  \texttt{http://hanson.gmu.edu}.

\item[] Hanson, R.\ (2006) ``Uncommon priors require origin disputes'',
  to appear in \textit{Theory and Decision}, available from 
  \texttt{http://hanson.gmu.edu}.

\item[] Harnad, S.\ (2001) ``What's wrong and right about Searle's Chinese 
  Room Argument?'' in: M.~Bishop and J.~Preston (editors) \textit{Essays on 
  Searle's Chinese Room Argument}, Oxford University Press.   Also
  available at \texttt{http://cogprints.org/4023/01/searlbook.htm}.

\item[] Knobe, J., Olum, K.~D., and Vilenkin, A.\ (2006) ``Philosophical
  implications of inflationary cosmology'', \textit{British Journal for the
  Philosophy of Science}, vol.~57, pp.~47-67.

\item[] Lewis, D.\ (1979) ``Attitudes \textit{de dicto} and \textit{de se}'',
  \textit{The Philosophical Review}, Vol.~88, pp.~513-543.

\item[] Lewis, D.\ (2001) ``Sleeping Beauty:\ reply to Elga'',
  \textit{Analysis}, vol.~61, pp.~171-176.

\item[] Leslie, J.\ (1996) \textit{The End of the World: The Science
  and Ethics of Human Extinction}, London: Routledge.

\item[] Marochnik, L.~S.\ (1983) ``On the origin of the solar system
  and the exceptional position of the sun in the galaxy'', \textit{Astrophysics
  and Space Science}, vol.~89, pp.~71-75.

\item[] Nozick (1969) ``Newcomb's problem and two principles of choice'',
  in N.~Rescher (editor) \textit{Essays in Honor of Carl G.~Hempel},
  Boston: Reidel, pp.~115-146.

\item[] Olum, K.~D.\ (2002) ``The doomsday argument and the number
  of possible observers'', \textit{The Philosophical Quarterly},
  vol.~52, pp.~164-184.

\item[] Olum, K.~D.\ (2004) ``Conflict between anthropic reasoning and
  observation'', \textit{Analysis}, vol.~64, pp.~1-8.

\item[] Smolin, L.\ (2006) ``Scientific alternatives to the anthropic
  principle'' (version 3, May 2006),  available from the e-print archive at 
  \texttt{http://arxiv.org/abs/hep-th/0407213} (first version
  posted July 2004).

\item[] Steinhardt, P.~J.\ and Turok, N.\ (2006) ``Why the cosmological
  constant is small and positive'', \textit{Science}, vol.~312, pp.~1180-1183.

\item[] Susskind, L.\ (2003) ``The anthropic landscape of string theory'',
  available from the e-print archive at
  {\texttt{http://arxiv.org/abs/hep-th/0302219}}.

\item[] Susskind, L.\ (2006) \textit{The Cosmic Landscape: String Theory
  and the Illusion of Intelligent Design}, New York: Little, Brown, 
  and Company.

\item[] Weinberg, S.\ (1987) ``Anthropic bound on the cosmological constant'',
  \textit{Physical Review Letters}, vol.~59, pp.~2607-2610.

\item[] Weinberg, S.\ (2000) ``\textit{A priori} probability distribution
  of the cosmological constant'', \textit{Physical Review D}, vol.~61,
  103505 (4 pages).

\item[] Whitman, W.~B., Coleman, D.~C., and Wiebe, W.~J.\ (1998)
  ``Prokaryotes: The unseen majority'', \textit{Proceedings of the
  National Academy of Sciences (USA)}, vol.~95, pp.~6578-6583.

\end{description}

\end{document}